\documentclass[reqno]{amsart}

\usepackage[dvipsnames]{xcolor}
\usepackage{amssymb,amscd,bbm,verbatim,graphicx,mathtools}
\usepackage[normalem]{ulem}
\usepackage[mathscr]{euscript}
\usepackage[colorlinks]{hyperref}
\newtheorem{thm}{Theorem}[section]

\newtheorem{lem}[thm]{Lemma}

\theoremstyle{definition}

\newtheorem{remark}[thm]{Remark}
\newtheorem{hyp}[thm]{Hypothesis}

\numberwithin{equation}{section}

\newcommand{\lVVert}{\left\vert\kern-0.25ex\left\vert\kern-0.25ex\left\vert}
\newcommand{\rVVert}{\right\vert\kern-0.25ex\right\vert\kern-0.25ex\right\vert}

\newcommand{\norm}[1]{\left\lVert #1\right\rVert}

\newcommand{\ov}{\overline}

\renewcommand{\Im}{\operatorname{Im}}
\newcommand{\supp}{\operatorname{supp}}

\newcommand{\iu}{\mathrm{i}}
\newcommand{\bm}{\textrm{m}}
\newcommand{\dom}{\operatorname{dom}}
\newcommand{\ran}{\operatorname{ran}}

\newcommand{\tr}{\operatorname{tr}}

\newcommand{\BVl}{\operatorname{BV}_{\rm loc}}

\newcommand{\sgn}{\operatorname{sgn}}
\newcommand{\diag}{\operatorname{diag}}
\newcommand{\sm}[1]{\big(\begin{smallmatrix}#1\end{smallmatrix}\big)}
\newcommand{\bb}[1]{\mathbb{#1}}
\newcommand{\mc}[1]{{\mathcal{#1}}}
\newcommand{\ms}[1]{{\mathscr{#1}}}

\newcommand{\id}{\mathbbm 1}
\newcommand{\pim}{\pi}
\newcommand{\Pim}{\Pi}

\newcommand{\ol}{{\overline{\lambda}}}
\newcommand{\la}{\lambda}
\newcommand{\loc}{{\rm loc}}
\newcommand{\mx}{{\rm max}}
\newcommand{\mn}{{\rm min}}

\newcommand{\<}{\langle}
\renewcommand{\>}{\rangle}

\newcommand{\disand}{\;\;\text{and}\;\;}

\begin{document}
\title[Fourier expansions]{On Fourier expansions for systems of ordinary differential equations with distributional coefficients}
\author{Steven Redolfi and Rudi Weikard}
\address{Department of Mathematics, University of Alabama at Birmingham, Birmingham, AL 35226-1170, USA}
\email{stevenre@uab.edu, weikard@uab.edu}
\date{\today}


\begin{abstract}
We study the spectral theory for the first-order system $Ju'+qu=wf$ of differential equations on the real interval $(a,b)$ where $J$ is a constant, invertible, skew-hermitian matrix and $q$ and $w$ are matrices whose entries are distributions of order $0$ with $q$ hermitian and $w$ non-negative.
Specifically, we construct a generalized Weyl-Titchmarsh $m$-function with corresponding spectral measure $\tau$ and a generalized Fourier transform after imposing certain conditions on $J$, $q$, and $w$.
Different conditions are motivated and studied in the later sections. A Fatou-type identity needed for our result is recorded in the appendix.
\end{abstract}
\maketitle

\section{Introduction}
In this paper we establish, under certain conditions, a Fourier expansion theorem for the differential equation
\begin{equation}\label{de}
Ju'+(q-\la w)u=wf
\end{equation}
posed on a real interval $(a,b)$.
Here $J$ is a constant, invertible, skew-hermitian $n\times n$-matrix and $\la$ is a complex parameter while $q$ and $w$ are $n\times n$-matrices whose entries are distributions of order $0$, $q$ is hermitian and $w$ non-negative.
If $f$ is in $L^2(w)$, a Hilbert space to be defined below, solutions $u$ of the equation are to be sought among functions of locally bounded variation.
We emphasize that, as the coefficients in equation \eqref{de} become rougher, so do the solutions $u$.
At some point it becomes impossible to define the products $qu$ and $wu$.
This point is reached when $q$ or $w$ are rougher than distributions of order $0$.

When the entries of $q$ and $w$ are locally integrable functions one has an existence and uniqueness theorem (the unique continuation property) for solutions of the initial value problem $Ju'+(q-\la w)u=wf$, $u(x_0)=u_0$.
This theorem is a central tool for the construction of Fourier expansions.
However, if $q$ or $w$ involve $\delta$-distributions, the existence or uniqueness of solutions of initial value problems can no longer be guaranteed.
In \cite{MR4047968} Ghatasheh and Weikard were able to overcome this obstacle in the case when the unique continuation property is only required for real $\la$.
It is the goal of this paper to weaken this requirement and thus to show how to establish a Fourier expansion under milder conditions.

The first to consider a Sturm-Liouville equation with measure coefficients was (to the best of our knowledge) Krein \cite{MR0054078} in 1952 when he modeled a vibrating string.
Further contributions we are aware of were made by Kac \cite{MR0080835}, Feller \cite{MR0068082}, Mingarelli \cite{MR706255}, Gesztesy and Holden \cite{MR914699}, Kurasov \cite{MR1397901}, Kurasov and Boman\cite{MR1443392}, Savchuk and Shkalikov \cite{MR1756602}, Eckhardt et al. \cite{MR3046408}, and Eckhardt and Teschl \cite{MR3095152}.
We refer the reader to Ghatasheh and Weikard \cite{MR4047968} and Eckhardt et al. \cite{MR3046408} for more details on the history of the subject.

We end this introduction with a short overview of the content of this paper.
Our main result is Theorem \ref{t:main} in Section \ref{S:FE} establishing the existence of the Fourier transform associated with a self-adjoint linear relation $T$ representing the differential equation \eqref{de} under appropriate boundary conditions.
The space of Fourier transforms of elements in $L^2(w)$ is a Hilbert space $L^2(\tau)$ where $\tau$ is a matrix-valued non-negative measure on $\bb R$.
This measure is determined through the Nevanlinna representation of a matrix-valued Weyl-Titchmarsh function $\bb M$ as we show in Section \ref{S:M}.
The matrix $\bb M$, in turn, is obtained by constructing Green's function, the kernel of the resolvent operator $R_\la=(T-\la)^{-1}$.
To this end we study the solutions of the equation $Ju'+(q-\la w)u=wf$ in Section \ref{S:sol} without the aid of the unique continuation property.
Before we can do all this we gather in Section \ref{S:pre} some material previously established in \cite{MR4047968}, \cite{MR4298818}, and \cite{MR4431055}.
In the appendix we record a result privately communicated to us by B. and C. Bennewitz \cite{BB-CB} concerning an extension of Fatou's theorem.

\section{Preliminaries}\label{S:pre}
In this section we gather some basic material on ordinary differential equations with distributional coefficients.

A distribution on the real interval $(a,b)$ is a linear functional $r$ on the set of infinitely often differentiable, compactly supported functions, the test functions, satisfying the following property: for any compact subinterval $K$ of $(a,b)$ there exist a positive constant $C$ and a non-negative integer $k$ such that
$$|r(\phi)|\leq C \sum_{j=0}^k \sup\{|\phi^{(j)}(x)|:x\in K\}$$
whenever $\phi$ is a test function with support in $K$.
If $k$ can be chosen to be $0$ independently of $K$ we say that $r$ is a distribution of order $0$.

Every distribution has a derivative and antiderivatives.
The derivative of $r$ is $\phi\mapsto r'(\phi)=-r(\phi')$.
According to du Bois-Reymond's lemma the difference of two antiderivatives of $r$ is a constant distributions, i.e., $\phi\mapsto C\int\phi\, dx$ for some $C\in\bb C$.
A distribution $r$ is called non-negative, if $r(\phi)\geq0$ whenever $\phi$ is a test function assuming only non-negative values.
The conjugate of a distribution $r$ is defined by $\phi\mapsto\ov{r}(\phi)=\ov{r(\ov\phi)}$.

Typical examples of distributions of order $0$ are (i) $\phi\mapsto \int_{(a,b)}\phi f\,dx$ when $f$ is a locally integrable function (with respect to Lebesgue measure) and (ii) $\phi\mapsto\phi(x_0)$, the $\delta$-distribution at $x_0$, when $x_0$ is some fixed point in $(a,b)$.
More generally, if $\mu$ is a non-negative Lebesgue-Stieltjes measure on $(a,b)$ and $h$ a function in $L^1_\loc(\mu)$, then $\phi\mapsto \int_{(a,b)} \phi h \mu$ is a distribution of order $0$.
Indeed, by Riesz's representation theorem, all distributions of order $0$ are of this type.
The cumulative distribution function of the (local) complex measure $h\mu$ is a function of locally bounded variation and the associated distribution is an antiderivative of $\phi\mapsto \int_{(a,b)} \phi h\mu$.
Conversely, given a function of locally bounded variation on $(a,b)$, it generates a Lebesgue-Stieltjes measure $\nu$ on each compact subset of $(a,b)$.
We may write $\nu=h|\nu|$ where $|\nu|$ is (locally) the total variation of $\nu$ and $h$ a function of modulus $1$ and thus in $L^1_\loc(|\nu|)$.
In the following we use the terms distribution of order $0$ and measure interchangeably thus deviating from the traditional meaning of the term measure.

Given a measure (or distribution of order $0$) $r$ and a function $f\in L^1_\loc(|r|)$ we define the distribution $rf=fr$ by setting $(rf)(\phi)=\int \phi f r$.
Note that $rf$ is also a distribution of order $0$.

Next we turn to matrices with distributional entries.
For any $n\times n$-matrix $w$ we have its adjoint $w^*$ defined in the canonical way as the conjugate of the transpose.
If $w=w^*$ then $w$ is called hermitian.
If, for any vector $z\in\bb C^n$, the distribution $z^*wz$ is non-negative, we call $w$ non-negative.
A non-negative distribution is necessarily hermitian.

Associated with a non-negative distribution is a non-negative measure and one may define the space $\mc L^2(w)$ consisting of those functions $f$  for which the measure $f^*wf$ is defined and finite.
In $\mc L^2(w)$ one may define the semi-scalar product $\<f,g\>=\int f^*wg$ and the semi-norm $\norm f=\<f,f\>^{1/2}$.
Identifying functions in $\mc L^2(w)$ whose difference have semi-norm $0$ we obtain the Hilbert space $L^2(w)$.

Throughout this paper we require the following canonical hypothesis to be valid.
More hypotheses (which we suspect to be only of a technical nature) need to be added later.

\begin{hyp}\label{hyp:m}
$(a,b)$ is a real interval.
$J$ is a constant, invertible, skew-hermitian $n\times n$-matrix.
$q$ and $w$ are $n\times n$-matrices whose entries are distributions of order $0$ on $(a,b)$; $q$ is hermitian and $w$ is non-negative.
\end{hyp}
Our goal is to investigate the differential equation
$$Ju'+(q-\la w)u=wf$$
where $\la\in\bb C$, $f\in L^2(w)$, and the solutions $u$ are to be sought among the balanced (to be defined presently) functions of locally bounded variation on $(a,b)$.
It makes sense to pose the equation, since all terms occurring in it are distributions of order~$0$.

For each function $u$ of locally bounded variation we denote by $u^-$ and $u^+$ its left- and right-hand limits, respectively.
More precisely, $u^-(x)=\lim_{t\uparrow x}u(t)$ and $u^+(x)=\lim_{t\downarrow x}u(t)$.
The balanced version $u^\#$ of a function $u$ of locally bounded variation is defined by $u^\#=(u^++u^-)/2$.
The reason for considering only balanced solutions of the differential equation lies in the integration by parts formula which states
\begin{equation}\label{ibypt}
\int_{[c,d]} (udv+vdu)=(uv)^+(d)-(uv)^-(c)+(2\kappa-1)\int_{[c,d]}(v^+-v^-)du
\end{equation}
whenever $[c,d]\subset(a,b)$, $u=\kappa u^++(1-\kappa)u^-$ and $v=\kappa v^++(1-\kappa)v^-$ for some fixed parameter $\kappa$.
Thus, only when $\kappa=1/2$, does the integration by parts formula have its familiar simple form when discontinuities of $u$ and $v$ are allowed.

If $r$ is a matrix-valued measure with cumulative distribution function $R$, let $\Delta_r(x)=r(\{x\})=R^+(x)-R^-(x)$.
If one of $\Delta_q(x)$ and $\Delta_w(w)$ is different from zero, satisfying the differential equation $Ju'+(q-\la w)u=wf$ requires that
$$B_+(x,\la)u^+(x)-B_-(x,\la)u^-(x)=\Delta_w(x)f(x)$$
where
$$B_\pm(x,\la)=J\pm\frac12(\Delta_q(x)-\la\Delta_w(x)).$$
If $B_+(x,\la)$ is not invertible it will not be possible to determine $u^+(x)$ uniquely even if $u^-(x)$ is given.
An analogous statement holds, of course, when $B_-(x,\la)$ is not invertible.
We emphasize that $B_-(x,\la)=-B_+(x,\ol)^*$.

Still we have the following theorem on existence and uniqueness of solutions of initial value problems for $Ju'+(q-\la w)u=wf$.

\begin{thm}[\cite{MR4047968}, Theorem 2.2]\label{EUIVP}
Suppose Hypothesis \ref{hyp:m} holds, that $f\in L^2(w)$ and $\la\in\bb C$, that the matrices $B_\pm(x,\la)$ are invertible for all $x\in(c,d)\subset(a,b)$, and that $x_0$ is a point in $(c,d)$.
Then the initial value problem $Ju'+(q-\la w)u=wf$, $u(x_0)=u_0\in\bb C^n$ has a unique balanced solution $u$ in the interval $(c,d)$.

We may pose an initial condition at $c$ (for $u^+$) or at $d$ (for $u^-$) if $q$ and $w$ are finite measures near the point in question.
\end{thm}

This theorem implies that we have fundamental matrices of solutions on $(c,d)$.
If $U$ is one such fundamental matrix of solutions for which $U^+$ or $U^-$ is equal to the identity matrix at some point $x_0\in(c,d)$ then it satisfies
the following Wronskian relationship
\begin{equation}\label{wronski}
U^\pm(x,\ol)^* J U^\pm(x,\la)=J \disand U^\pm(x,\la) J^{-1} U^\pm(x,\ol)^*=J^{-1}
\end{equation}
which follows from Lemma 3.2 in \cite{MR4047968}.

Let us quickly recall from \cite{MR4431055} some basic facts about the matrices $B_\pm(x,\la)$.
Let $\Lambda_x$ denote the set of those $\lambda\in\bb C$ where invertibility of either $B_+(x,\la)$ or $B_-(x,\la)$ fails.
The set of those $x\in(a,b)$ where invertibility of $B_+(x,\la)$ or $B_-(x,\la)$ fails is denoted by $\Xi_\lambda$.
Then $\Lambda_x$ is empty for all but countably many $x\in(a,b)$.
When it is not empty, it is either finite or all of $\bb C$.
Moreover, it is symmetric with respect to the real axis.
The intersections of $\Xi_\lambda=\Xi_\ol$ with any compact subset of $(a,b)$ are finite.

Associated with the differential equation $Ju'+qu=wf$ are the linear relations
$$\mc T_\mx=\{(u,f)\in \mc L^2(w)\times \mc L^2(w): u\in\BVl^\#((a,b))^n, Ju'+qu=wf\}$$
and
$$\mc T_\mn=\{(u,f)\in \mc T_\mx: \text{$\supp u$ is compact in $(a,b)$}\}.$$
Here $\BVl^\#((a,b))^n$ denotes the space of $\bb C^n$-valued functions on $(a,b)$ which are of locally bounded variation and balanced.
The relations $\mc T_\mx$ and $\mc T_\mn$ give rise to linear relations in the Hilbert space $L^2(w)\times L^2(w)$
$$T_\mx=\{([u],[f])\in L^2(w)\times L^2(w):(u,f)\in\mc T_\mx\}$$
and
$$T_\mn=\{([u],[f])\in L^2(w)\times L^2(w):(u,f)\in\mc T_\mn\}.$$
$T_\mn$ and $T_\mx$ are called minimal and maximal relation, respectively.

In \cite{MR4298818} it was shown that $T_\mn^*=T_\mx$ proving that $T_\mn$ is a symmetric relation.
Therefore we have von Neumann's relation
$$T_\mx=\ov{T_\mn}\oplus D_\iu\oplus D_{-\iu}$$
where $\ov T_\mn$ is the closure of $T_\mn$ and $D_\la=\{(u,\la u)\in T_\mx\}$.
The dimensions of $D_\la$ are constant as $\la$ varies in either the lower or upper half plane.
These numbers are called deficiency indices of $T_\mn$ and we use the notation $n_\pm=\dim D_{\pm\iu}$.
While, in general, the equation $Ju'+qu=\la wu$ may have infinitely many linearly independent solutions we showed in \cite{MR4431055} that the deficiency indices are still no larger than $n$.
We also established Green's formula or Lagrange's identity, i.e.,
\begin{equation*}\label{Lagrange}
(v^*Ju)^-(b)-(v^*Ju)^+(a)=\<v,f\>-\<g,u\>
\end{equation*}
provided that $(v,g)$ and $(u,f)$ are in $\mc T_\mx$.

We showed in \cite{MR4431055} that the closed symmetric extensions of $T_\mn$ are given as restrictions of $T_\mx$ by boundary conditions.
In particular, this is true for the self-adjoint extensions in which case we have the following special case of Theorem 4.4 in \cite{MR4431055}.
\begin{thm}\label{T:3.5}
If the deficiency indices of $T_\mn$ satisfy $n_+=n_-$ and if $(v_1,g_1)$, ..., $(v_{n_\pm},g_{n_\pm})$ are linearly independent elements of $D_\iu\oplus D_{-\iu}$ such that
\begin{equation}\label{200820.1}
(g_k^*Jg_\ell)^-(b)-(g_k^*Jg_\ell)^+(a)=0\;\; \text{for $1\leq k,\ell\leq n_\pm$,}
\end{equation}
then
\begin{equation}\label{200821.1}
T=\{(u,f)\in T_\mx: (g_j^*Ju)^-(b)-(g_j^*Ju)^+(a)=0\text{ for $j=1,...,n_\pm$}\}
\end{equation}
is a self-adjoint extension of $T_\mn$.

Conversely, if $T$ is a self-adjoint extension of $T_\mn$, then $T$ is given by \eqref{200821.1} for appropriate elements $(v_1,g_1)$, ..., $(v_{n_\pm},g_{n_\pm})$ of $D_\iu\oplus D_{-\iu}$ for which \eqref{200820.1} holds.
\end{thm}

If $T$ is a self-adjoint relation and $\la$ is in $\varrho(T)$, the resolvent set of $T$, we showed in \cite{MR4431055} that the resolvent operator $R_\la=(T-\la)^{-1}$ is an integral operator, i.e., we established the existence of a Green's function for $T$.

Some solutions of $Ju'+qu=\la w u$ may have norm zero, i.e., $\norm{u}^2=\int u^*wu=0$.
Since $w$ is non-negative $\norm u=0$ if and only if the distribution $wu=0$.
It follows that such a function $u$ satisfies $Ju'+qu=\la w u$ for any $\lambda$, in particular for $\lambda=0$.
The space of solutions of $Ju'+qu=0$ satisfying $wu=0$ is denoted by $\mc L_0$.

In particular, if $u\in \mc L_0$ (and hence $[u]=0$ and $f=0$), then
\begin{equation}\label{lagrange&L0}
(v^*Ju)^-(b)-(v^*Ju)^+(a)=\<v,f\>-\<g,u\>=0
\end{equation}
for all $(v,g)\in\mc T_\mx$.
Since $(v,g)\in D_\iu\oplus D_{-\iu}$ implies that $(g,-v)\in D_\iu\oplus D_{-\iu}$, it follows that any element of $\mc L_0$ satisfies the boundary conditions defining a self-adjoint relation $T$.

\section{Constructing solutions}\label{S:sol}
Throughout the remainder of the paper we require, in addition to Hypothesis~\ref{hyp:m}, that $\Lambda_x\cap\bb R=\emptyset$ for all but finitely many $x\in(a,b)$.
The points where $\Lambda_x\cap\bb R\neq\emptyset$ are (included among) the points $x_1$, ..., $x_N$ which we consider ordered by size.
We shall allow that $N=0$ when $\Lambda_x\cap\bb R$ is empty for all $x\in(a,b)$.
We also define $x_0=a$ and $x_{N+1}=b$.
Note that, while the set $\Lambda=\bigcup_{x\in(a,b)}\Lambda_{x}$ may be all of $\bb C$, the set $$\tilde\Lambda=\bigcup_{x\not\in\{x_1,...,x_N\}}\Lambda_{x}$$
does not intersect $\bb R$ and is countable.

According to Theorem \ref{EUIVP} initial value problems for $Ju'+qu=\la wu$ have unique solutions on $(x_j,x_{j+1})$ for $j=0,..., N$.
In each interval $(x_j,x_{j+1})$ we choose a balanced fundamental system of solutions $U_j(\cdot,\la)$ such that $U_j(\xi_j,\la)=\id$ for some point $\xi_j$ which is a point of continuity for $q$ and $w$, i.e., $q(\{\xi_j\})=w(\{\xi_j\})=0$.
These are extended to all of $(a,b)$ as balanced functions which vanish outside $\ov{(x_j,x_{j+1})}$.
The extensions are also denoted by $U_j(\cdot,\la)$.
With these we define the $n\times n(N+1)$-matrix.
$$\ms U(x,\la)=(U_0(x,\la), ..., U_N(x,\la)).$$
We may now define the Fourier transform $\ms Ff$ of an element $f\in L^2(w)$, at least as long as $f$ is compactly supported, by
$$(\ms Ff)(\la)=\int \ms U(\cdot,\ol)^*wf,$$
a vector in $\bb C^{n(N+1)}$.

At this point let us introduce some more notation needed below.
Identity matrices of various sizes are denoted by $\id$ while (square or rectangular) zero matrices are denoted by $0$.
Sometimes it may be advisable to indicate dimensions as subscripts, e.g., $\id_k$ is the $k\times k$ identity matrix while $0_{k,j}$ the $k\times j$ zero matrix.

In \cite{MR4298818} and \cite{MR4431055} the equation  $Ju'+(q-\la w)u=wf$ was studied in the presence of points across which solutions cannot be uniquely continued.
There, the description of solutions involved a two-diagonal block matrix $\bb B$ with blocks of size $n\times n$.
Such a two-diagonal block matrix has the form $A\bb E_\top+B\bb E_\bot$ where $A$ and $B$ are $N\times N$ diagonal block matrices while  $\bb E_\top$ and $\bb E_\bot$ are $N\times (N+1)$-block matrices, which, respectively, strip the first and last $n$ components off a vector in their domain $\bb C^{n(N+1)}$.
Thus $\bb E_\top=(0_{nN,n},\id_{nN})$ and $\bb E_\bot=(\id_{nN},0_{nN,n})$.
In our particular application we have
\[\bb B(\lambda)=\mc B(\la)\mc U^+(\lambda)\bb E_\top+\mc B(\ol)^*\mc U^-(\lambda)\bb E_\bot\]
where
\[\mc B(\lambda)=\diag(B_+(x_1,\lambda),...,B_+(x_N,\lambda)),\]
\[\mc U^-(\lambda)=\diag(U_0^-(x_1,\lambda),...,U_{N-1}^-(x_N,\lambda)),\]
and
\[\mc U^+(\lambda)=\diag(U_1^+(x_1,\lambda),...,U_{N}^+(x_N,\lambda)).\]

We will also need the matrix
\[\tilde{\bb B}(\la)=\mc B(\la)\mc U^+(\lambda)\bb E_\top-\mc B(\ol)^*\mc U^-(\lambda)\bb E_\bot.\]
If $N=0$, $\bb E_\top$ and $\bb E_\bot$ have to be considered as linear transformations from $\bb C^n$ to $\{0\}$.
In this case $\bb B(\la)$ and $\tilde{\bb B}(\la)$ are ``matrices'' with no rows.
This causes no problems below.
We also let $\mc J=\diag(J,..., J)$ stand for either an $N\times N$ or an $(N+1)\times (N+1)$ block matrix.

\subsection{A representative of \texorpdfstring{$R_\la f$}{Rf}}
If $\la\in\varrho(T)\setminus\tilde\Lambda$ and $f\in L^2(w)$ is compactly supported, we shall now construct a balanced representative of $R_\la f$.

\subsubsection{General solution}
Using the variation of constants formula in the intervals $(x_{j-1},x_j)$ shows that any solution $u$ of $Ju'+(q-\la w)u=wf$ satisfies
\begin{equation}\label{201106.2}
u^-(x)=\ms U^-(x,\la)\big(\tilde u + \mc J^{-1}\int_{(a,x)}\ms U(\cdot,\ol)^*wf\big)
\end{equation}
where $\tilde u=(\tilde u_0, ...,\tilde u_N)^\diamond$ is an appropriate vector in $\bb C^{n(N+1)}$.%
\footnote{We use ${}^\diamond$ to indicate a block column vector as in $(c_1, ..., c_N)^\diamond=(c_1^\top, ..., c_N^\top)^\top$.}
For $u$ to be indeed a solution it is necessary and sufficient that
\begin{equation}\label{220521.1}
-B_-(x_j,\la)u^-(x_j)+B_+(x_j,\la)u^+(x_j)=\Delta_w(x_j) f(x_j)
\end{equation}
for $j=1,...,N$.
Since, by equation \eqref{wronski}, $U_j^-(x_j,\la)J^{-1}U_j^-(x_j,\ol)^*=J^{-1}$ we find
\begin{multline*}
u^-(x_j)=U_{j-1}^-(x_j,\la)\big(\tilde u_{j-1}+J^{-1}\int_{(a,x_j)}U_{j-1}(\cdot,\ol)^*wf\big)\\
 =U_{j-1}^-(x_j,\la)\big(\tilde u_{j-1}+J^{-1}\int U_{j-1}(\cdot,\ol)^*wf\big)-\frac12 J^{-1}(\Delta_w f)(x_j)\\
\end{multline*}
so that $(u^-(x_1),..., u^-(x_N))^\diamond=\mc U^-(\la)\bb E_\bot\big(\tilde u+\mc J^{-1}(\ms Ff)(\la) \big)-\frac12 \mc J^{-1}\ms W_0(f)$
where $\ms W_0(f)=((\Delta_wf)(x_1),...,(\Delta_wf)(x_N))^\diamond$.

Similarly, $(u^+(x_1),..., u^+(x_N))^\diamond=\mc U^+(\la)\bb E_\top\tilde u+\frac12\mc J^{-1}\ms W_0(f)$, since
\begin{multline*}
u^+(x_j)=U_{j}^+(x_j,\la)\big(\tilde u_{j}+J^{-1}\int_{(a,x_j]}U_{j}(\cdot,\ol)^*wf\big)\\
 =U_{j}^+(x_j,\la)\tilde u_{j}+\frac12 J^{-1}(\Delta_wf)(x_j).
\end{multline*}
Thus, taking into account that $\mc B(\la)-\mc B(\ol)^*=2\mc J$, equations \eqref{220521.1} show that $\tilde u$ must solve the system
\begin{equation}\label{eq:bb}
\bb B(\la)\tilde u=-\mc B(\ol)^*\mc U^-(\la)\bb E_\bot \mc J^{-1}(\ms Ff)(\la).
\end{equation}

\subsubsection{Integrability conditions}
Next we require that $u$ is an element of $\mc L^2(w)$.
We define $N_+(\la)\subset\bb C^n$ to be the set of all $\eta$ such that $u=U_N(\cdot,\la)\eta$ satisfies $\int_{(x_N,b)}u^*wu<\infty$.
Similarly, $N_-(\la)\subset\bb C^n$ is the set of all $\eta$ such that $u=U_0(\cdot,\la)\eta$ satisfies $\int_{(a,x_1)}u^*wu<\infty$.
We also denote the orthogonal projections onto $N_\pm(\la)$ by $P_\pm(\la)$, respectively.

To pick out the solutions of the differential equation which are in $L^2(w)$ we proceed as follows.
Recall that $wf$ is compactly supported.
In $(a,x_1)$ we have
$$u^-(x)=U_0^-(x,\la)(\tilde u_0+J^{-1}\int_{(a,x)}U_0(\cdot,\ol)^*wf).$$
In particular, if $x$ is below the support of $wf$, we have $u^-(x)=U_0^-(x,\la)\tilde u_0$.
Hence we want $\tilde u_0\in N_-(\la)$.
Define $\ms Q_-(\la)=(\id-P_-(\la),0,...,0)\in\bb C^{n\times n(N+1)}$ so that the condition becomes
\begin{equation}\label{eq:pp}
\ms Q_-(\la)\tilde u=0.
\end{equation}
For $x\in(x_N,b)$ we get instead
$$u^-(x)=U_N^-(x,\la)(\tilde u_N+J^{-1}\int_{(x_N,x)}U_N(\cdot,\ol)^*wf)$$
In particular, $u^-(x)=U_N^-(x,\la){(\tilde u_N+J^{-1}(\ms Ff)_N(\la))}$ if $x$ is above the support of $wf$.
Hence $(\id-P_+(\la))\tilde u_N=-(\id-P_+(\la))J^{-1}(\ms Ff)_N(\la)$.
Setting $\ms Q_+(\la)=(0,...,0,\id-P_+(\la))\in\bb C^{n\times n(N+1)}$ the condition becomes
\begin{equation}\label{eq:pm}
\ms Q_+(\la)\tilde u=-\ms Q_+(\la)\mc J^{-1}(\ms Ff)(\la).
\end{equation}

If $a$ (or $b$) is regular, then $P_-(\la)$ (or $P_+(\la)$) is the identity and the corresponding condition is vacuous.

\subsubsection{Boundary conditions}
We now invoke Theorem \ref{T:3.5} to deal with the boundary conditions representatives of $R_\la$ have to satisfy.
These boundary conditions are given by
\begin{equation}\label{sabc}
(g_j^*Ju)^-(b)-(g_j^*Ju)^+(a)=0, \;\; \text{for $j=1,...,n_\pm$}
\end{equation}
where $(v_1,g_1)$, ..., $(v_{n_\pm},g_{n_\pm})$ are linearly independent elements of $D_\iu\oplus D_{-\iu}$ satisfying equation \eqref{200820.1}.
Let $g=(g_1,...,g_{n_+})$ and introduce the matrices\footnote{The definition of $A_-$ is differs from the one in \cite{MR4047968} by a factor of $-1$.}
$$A_-(\la)=-(g^*JU_0(\cdot,\la)P_-(\la))^+(a),$$
$$A_+(\la)=(g^*JU_N(\cdot,\la)P_+(\la))^-(b),$$
and
$$\ms A_-(\la)=(A_-(\la),0,...,0),\quad \ms A_+(\la)=(0,...,0,A_+(\la))$$
with $N$ blocks of zero-matrices.
Then the boundary conditions are
\begin{equation}\label{eq:bc}
(\ms A_+(\la)+\ms A_-(\la))\tilde u=-\ms A_+(\la)\mc J^{-1}(\ms Ff)(\la).
\end{equation}

\subsubsection{Solutions of zero norm}\label{S:3.4}
Finally, we take account of the presence of solutions with norm $0$.
We define $N_0=\{\eta\in\bb C^{n(N+1)}: \ms U(\cdot,0)\eta\in\mc L_0\}$.
Adding any element of $\mc L_0$ to a solution $u$ of $Ju'+(q-\la w)u=wf$ which is in $\mc L^2(w)$ and satisfies the boundary conditions will yield a solution with the same properties.
Each of these is given by \eqref{201106.2} for an appropriate choice of $\tilde u$ and exactly one of them will have $\tilde u\in N_0^\perp$.
Thus, if $\bb P$ is the orthogonal projection onto $N_0^\perp$, we require
\begin{equation}\label{eq:ppp}
(\id-\bb P)\tilde u=0.
\end{equation}
Note that $\bb P$ does not depend on $\la$.

We also emphasize that $w\ms U(\cdot,\la)(\id-\bb P)\eta$ is the zero distribution for any $\eta\in\bb C^{n(N+1)}$ and any $\la\in\bb C\setminus\tilde\Lambda$.
Thus $w\ms U(\cdot,\la)=w\ms U(\cdot,\la)\bb P$ or, equivalently, $\ms U(\cdot,\ol)^*w=\bb P\ms U(\cdot,\ol)^*w$ and, in particular,
\begin{equation}\label{puw}
(\ms Ff)(\la)=\bb P(\ms Ff)(\la)
\end{equation}
which we record here for later use.

We shall also need the following.
\begin{lem}\label{L:3.1}
$\bb B(\la)(\id-\bb P) = 0$, $\ms Q_\pm(\la)(\id-\bb P) = 0$, and $(\ms A_+(\la) + \ms A_-(\la))(\id-\bb P) = 0$.
\end{lem}

\begin{proof}
Suppose $v$ is an arbitrary element of $\bb C^{n(N+1)}$.
Then $(\id-\bb P)v$ is in $N_0$ and $\ms U(\cdot,\la)(\id-\bb P)v$ is an element of $\mc L_0$.
As such it solves the differential equation $Ju'+qu=\la wu$ which implies that $\bb B(\la)(\id-\bb P)v=0$ according to equation \eqref{eq:bb} for $f=0$.
It is also (trivially) in $L^2(w)$ and satisfies the boundary conditions, as we argued in equation \eqref{lagrange&L0}, implying the other two claims. \end{proof}

\subsubsection{Putting it all together}
We now collect the equations \eqref{eq:bb}, \eqref{eq:pp}, \eqref{eq:pm}, \eqref{eq:bc}, and \eqref{eq:ppp} into the system
\begin{equation}\label{201106.1}
\bb F(\la)\tilde u=\bb H_\ell(\la)\mc J^{-1}(\ms Ff)(\la)
\end{equation}
where
$$\bb F(\la)=\begin{pmatrix}
  \bb B(\la)\\
  \ms Q_-(\la)\\
  \ms Q_+(\la)\\
  \ms A_+(\la)+\ms A_-(\la)\\
  \id-\bb P
\end{pmatrix}
\disand
\bb H_\ell(\la)=-\begin{pmatrix}
  \mc B(\ol)^*\mc U^-(\lambda)\bb E_\bot\\
  0\\
  \ms Q_+(\la)\\
  \ms A_+(\la)\\
  0
\end{pmatrix}.$$

We will show presently that $\bb F$ has trivial kernel and hence full column rank.
Therefore it has a left inverse $\bb F^\dag$ and since we know that a solution exists we may solve \eqref{201106.1} for $\tilde u$.
Using $\tilde u=\bb P\tilde u$ and equation \eqref{puw} this gives
\begin{equation}\label{eq:tu}
\tilde u=\bb P\bb F(\la)^\dag \bb H_\ell(\la)\mc J^{-1}\bb P(\ms Ff)(\la)=\bb M_\ell(\la)(\ms Ff)(\la)
\end{equation}
thereby defining the matrix $\bb M_\ell$.
With this value of $\tilde u$ we repeat equation \eqref{201106.2} to get
\begin{equation}\label{finum}
u^-(x)=\ms U^-(x,\la)\int (\bb M_\ell(\la)+\mc J^{-1}\chi_{(a,x)})\ms U(\cdot,\ol)^*wf.
\end{equation}

\begin{lem}
If $\lambda\in\bb C\setminus\tilde{\Lambda}$ and $\Im(\la)\neq0$, then $\ker\mathbb{F}(\lambda)=\{0\}$.
\end{lem}

\begin{proof}
Assume $v\in\ker\mathbb{F}(\lambda)$.
Since then $v\in\ker\mathbb{B}(\lambda)$ it follows that $u=\ms{U}(\cdot,\lambda)v$ is a solution of $Ju'+qu=\lambda wu$.
Notice further that $u$ must be in $\mathcal{L}^2(w)$ since $\ms{Q}_\pm(\lambda)v=0$.
Finally, $(\ms{A}_+(\lambda)+\ms{A}_-(\lambda))v=0$ gives that $u$ satisfies the boundary conditions.
It is thus the case that $([u],\lambda [u])\in T$.
Since $\lambda$ is not real we must have $[u]=0$ so that $u\in\mathcal{L}_0$.
We thus see that $v\in N_0$.
However, we also have $(\id-\mathbb{P}) v=0$ and hence $v=0$.
\end{proof}

We now obtain $u^+$ by taking limits from the right in equation \eqref{finum} to get
$$u^+(x)=\ms U^+(x,\la)\int (\bb M_\ell(\la)+\mc J^{-1}\chi_{(a,x]})\ms U(\cdot,\ol)^*wf.$$
Hence
\begin{multline}
u(x)=\ms U(x,\la)\int (\bb M_\ell(\la)+\mc J^{-1}\chi_{(a,x)})\ms U(\cdot,\ol)^*wf\\
 +\frac12\ms U^+(x,\la)\mc J^{-1}\ms U(x,\ol)^*\Delta_w(x)f(x)\label{eq:u}
\end{multline}
is a representative of $R_\la f$.

\subsection{Striving for symmetry}
We now construct another representative $v$ of $R_\la f$ by starting from integrals over $(x,b)$ instead of $(a,x)$.
Specifically,
$$v^+(x)=\ms U^+(x,\la)\big(\tilde v - \mc J^{-1}\int_{(x,b)}\ms U(\cdot,\ol)^*wf\big)$$
with an appropriate vector $\tilde v$ in $\bb C^{n(N+1)}$.

Equation \eqref{201106.1} becomes
\begin{equation}\label{221107.1}
\bb F(\la)\tilde v=\bb H_r(\la)\mc J^{-1}(\ms Ff)(\la)
\end{equation}
where $\bb F$ is as before but
$$\bb H_r(\la)=\begin{pmatrix}
  \mc B(\la)\mc U^+(\lambda)\bb E_\top\\
  \ms Q_-(\la)\\
  0\\
  \ms A_-(\la)\\
  0
\end{pmatrix}.$$
Analogously to \eqref{eq:tu} we get now
\begin{equation}\label{eq:tv}
\tilde v=\bb P\bb F(\la)^\dag \bb H_r(\la)\mc J^{-1}\bb P(\ms Ff)(\la)=\bb M_r(\la)(\ms Ff)(\la)
\end{equation}
defining $\bb M_r(\la)$.
Taking limits of $v^+$ from the left we obtain $v^-$ and therefore, in analogy to \eqref{eq:u},
\begin{multline}\label{eq:v}
v(x)=\ms U(x,\la)\int (\bb M_r(\la)-\mc J^{-1}\chi_{(x,b)})\ms U(\cdot,\ol)^*wf\\
 -\frac12\ms U^-(x,\la)\mc J^{-1}\ms U(x,\ol)^*\Delta_w(x)f(x).
\end{multline}

We emphasize that, since $u$ and $v$ are representatives of $R_\la f$, their difference must be an element of $\mc L_0$.

\subsection{Constructing \texorpdfstring{$\bb M$}{M}}
We now define $\bb M(\la)=\frac12(\bb M_\ell(\la)+\bb M_r(\la))$, i.e.,
$$\bb M(\la)=\bb P\bb F(\la)^\dag\bb H(\la)\mc J^{-1}\bb P$$
where
$$\bb F(\la)=\begin{pmatrix}
  \bb B(\la)\\
  \ms Q_-(\la)\\
  \ms Q_+(\la)\\
  \ms A_+(\la)+\ms A_-(\la)\\
  \id-\bb P
\end{pmatrix},
\quad
\bb H(\la)=\frac{\bb H_\ell(\la)+\bb H_r(\la)}2=
 \frac12\begin{pmatrix}
  \tilde{\bb B}(\la)\\
  \ms Q_-(\la)\\
 -\ms Q_+(\la)\\
  \ms A_-(\la)-\ms A_+(\la)\\
  0
\end{pmatrix}.$$

Then we find, combining \eqref{eq:u} and \eqref{eq:v}, that
\begin{multline}
(\mc R_\la f)(x)=\frac{u(x)+v(x)}2
 =\ms U(x,\la)\int (\bb M(\la)+\frac12\mc J^{-1}\sgn(x-\cdot))\ms U(\cdot,\ol)^*wf\\
     +\frac14 (\ms U^+(x,\la)-\ms U^-(x,\la))\mc J^{-1}\ms U(x,\ol)^*\Delta_w(x)f(x).
     \label{elf}
\end{multline}
This is our final representative of $R_\la f$.
We emphasize that the last term in \eqref{elf} is zero for all but countably many $x\in(a,b)$.

\section{Properties of \texorpdfstring{$\bb M$}{M}}\label{S:M}
The function
$$\bb M(\la)=\bb P\bb F(\la)^\dag\bb H(\la)\mc J^{-1}\bb P,$$
occurring in equation \eqref{elf}, is a matrix-valued Nevanlinna function in a wide variety of circumstances (including all regular cases and when $N=0$) as we will show in this section.
In fact we believe this is always the case but can, at this point, not prove it.

Recall that a Nevanlinna function $\bb M$ is a function defined on $\bb C\setminus\bb R$ with the following properties: (i) $\bb M$ is symmetric in the sense that $\bb M(\la)=\bb M(\ol)^*$, (ii) $\Im \bb M$ is non-negative in the upper half-plane, and (iii) $\bb M$ is analytic.
As a Nevanlinna function $\bb M$ has the representation
\begin{equation}\label{mnl}
\bb M(\la)=A+B\la+\int_{\bb R}\Big(\frac{1}{t-\la}-\frac{t}{t^2+1}\Big)\tau(t)
\end{equation}
where $A$ is hermitian, $B$ is non-negative, and $\tau$ is a non-negative matrix-valued measure whose significance will become clear in the next section.
The measure $\tau$ is called a spectral measure.

\subsection{Symmetry of \texorpdfstring{$\bb M$}{M}}
For $\la\in\bb C\setminus\tilde\Lambda$ define the set
$$\mathbf{B}(\la)=\{(\ms Ff)(\la): \text{$f\in L^2(w)$, $\supp f$ compact}\}.$$
\begin{lem}\label{L:4.1}
If $\la,\mu\in\bb C\setminus\tilde\Lambda$ then $\mathbf{B}(\la)=\mathbf{B}(\mu)$.
\end{lem}
\begin{proof}
Let $u=\ms U(\cdot,\ol)\alpha$ and $v=\ms U(\cdot,\overline\mu)\alpha$.
Since $\alpha\in \mathbf{B}(\la)^\perp$ if and only if $wu=0$ it follows for such $\alpha$ that $u$ satisfies the differential equation $Ju'+qu=\ol wu=0=\overline\mu wu$ on $(x_j,x_{j+1})$.
Since $u(\xi_j)=v(\xi_j)=\alpha_j$ we have actually that $u$ and $v$ coincide on $(x_j,x_{j+1})$ and hence everywhere.
Therefore $wv=0$ and $\alpha\in \mathbf{B}(\mu)^\perp$.
Switching the roles played by $\la$ and $\mu$ shows the other inclusion.
\end{proof}
In view of Lemma \ref{L:4.1} we will write subsequently simply $\mathbf B$ rather than $\mathbf{B}(\la)$.

Because of equation \eqref{puw} it is clear that $\mathbf{B}$ is a subspace of $N_0^\perp=\ran \bb P$.
This inclusion may be strict\footnote{Let, for instance, $(a,b)=\bb R$, $J=\sm{0&-1\\ 1&0}$, $q=\sm{0&0\\ 0&2}\delta_0$, and $w=\sm{2&0\\ 0&0}\delta_0$.
Then $\dim\mathbf{B}=1$ while $\dim\ran\bb P=3$.}
but when it is not we have the following lemma.

\begin{lem}\label{L.4.1}
Suppose $\la\in\varrho(T)\setminus\tilde\Lambda$ and $\mathbf{B}=\ran \bb P$.
Then $\bb M(\ol)^*=\bb M(\la)$.
This is true, in particular, when $N=0$.
\end{lem}

\begin{proof}
Since $R_\ol=R_\la^*$ we have
$$0=\<g,R_\la f\>-\<R_\ol g,f\>
=(\ms Fg)(\ol)^*(\bb M(\la)-\bb M(\ol)^*)(\ms Ff)(\la).$$
Here we used \eqref{wronski} to show that the integrated terms stemming from \eqref{elf} cancel each other.
The claim follows then because $(\ran\bb P)^\perp=N_0$ is in the kernel of both $\bb M(\la)$ and $\bb M(\ol)^*$.

Next we show that for $N=0$ we have indeed $\mathbf{B}=\ran \bb P$.
Let $\alpha\in \mathbf{B}^\perp$ so that $\int f^*w \ms U(\cdot,\ol)\alpha=0$.
But here $\ms U(\cdot,\ol)=U_0(\cdot,\ol)$ implying that $\ms U(\cdot,\ol)\alpha$ is a solution of $Ju'+qu=\ol wu$ which is perpendicular to a dense set in $L^2(w)$.
This shows that $\alpha\in N_0=(\ran\bb P)^\perp$.
\end{proof}

If $N\geq 1$ we still have that $u=\ms U(\cdot,\ol)\alpha$ has norm $0$ when $\alpha\in \mathbf{B}^\perp$ but it may not be a solution of the differential equation anymore.
To cover the exceptional cases where $\mathbf{B}\neq\ran \bb P$ we introduce the matrix
$$\Omega(\la)=\bb H(\la)\mc J^{-1}\bb P\bb F(\ol)^*+\bb F(\la)\bb P\mc J^{-1}\bb H(\ol)^*.$$
Applying $\bb P \bb F(\la)^\dag$ to $\Omega(\la)$ from the left and $\bb F(\ol)^{\dag*}\bb P$ from the right shows that $\Omega(\la)=0$ implies $\bb M(\la)=\bb{M}(\ol)^*$, a result we state as a lemma.

\begin{lem}
If $\la\in\varrho(T)\setminus\tilde\Lambda$ and $\Omega(\la)=0$, then $\bb M(\la)=\bb M(\ol)^*$.
\end{lem}

In general $\Omega$ is a $5\times 5$-block-matrix and we mean rows and columns of blocks when we speak simply of rows and columns in the following.
Note that for $N=0$ the first row and column are actually absent, while for $n_\pm=0$ boundary conditions are not needed so that the fourth row and column are absent.
For simplicity of notation we will, nevertheless, maintain the labeling of the remaining blocks as if those others were present.

Since the last rows of blocks in both $\bb H$ and $\bb F\bb P$ are zero it follows that the last row of $\Omega$ is also zero.
Since $\Omega(\ol)^*=-\Omega(\la)$ the same is true for the last column of $\Omega$.
To investigate the remaining $4\times 4$-block-matrix we introduce
$$X(\la)=\begin{pmatrix}
  \mc B(\la)\mc U^+(\lambda)\bb E_\top\\
  \ms Q_-(\la)\\
  0\\
  \ms A_-(\la)\\
\end{pmatrix}
\disand
Y(\la)=\begin{pmatrix}
  \mc B(\ol)^*\mc U^-(\lambda)\bb E_\bot\\
  0\\
  \ms Q_+(\la)\\
  \ms A_+(\la)\\
\end{pmatrix}.$$
From Lemma \ref{L:3.1} we obtain that $(X + Y)\bb P = X + Y$.
Thus, taking cancellations into account, we find that the block $\Omega_{\ell,k}$ is given by
$$\Omega_{\ell,k}(\la)=X_\ell(\la)\mc J^{-1}X_k(\ol)^*-Y_\ell(\la)\mc J^{-1}Y_k(\ol)^*$$
for $1\leq\ell,k\leq4$.

\begin{lem}
If $\la\in\varrho(T)\setminus\tilde\Lambda$ the first row and column of $\Omega(\la)$ (if present) are zero.

\end{lem}

\begin{proof}
$\Omega_{1,1}(\la)$ involves the term
$$\mc U^+(\la)\bb{E}_\top\mathcal{J}^{-1}\bb{E}_{\top}^*\mc U^+(\ol)^*$$
where $\mc J^{-1}$ is an $n(N+1)\times n(N+1)$ matrix.
In view of equation \eqref{wronski} and the structure of $\bb E_\top$ this product is equal to the $nN\times nN$ matrix $\mc J^{-1}$.
Using the analogous argument for the term involving $\bb E_\bot$ we obtain
$$\Omega_{1,1}(\la)=\mathcal{B}(\lambda)\mathcal{J}^{-1}\mathcal{B}(\ol)^*-\mathcal{B}(\ol)^*\mathcal{J}^{-1}\mathcal{B}(\lambda),$$
a block-diagonal matrix.
Since $B_+(\cdot,\ol)^*=-B_-(\cdot,\lambda)$ and $B_\pm=J\pm A$ for a suitable matrix $A$ we get $B_-(x_k,\lambda)J^{-1}B_+(x_k,\lambda)-B_+(x_k,\lambda)J^{-1}B_-(x_k,\lambda)=0$ and hence $\Omega_{1,1}(\la)=0$.

For $\ell=2$ and $k=1$ we find
$$\Omega_{2,1}(\la)=\ms Q_-(\la)\mathcal{J}^{-1}\mathbb{E}_\top^*\mathcal{U}^+(\ol)^*\mathcal{B}(\ol)^*.$$
This is $0$ since only the first block in the row vector $\ms Q_-(\la)\mathcal{J}^{-1}$ is non-zero and therefore annihilated by $\mathbb{E}_\top^*$.
A similar argument works to show that $\Omega_{3,1}(\la)=0$ and $\Omega_{4,1}(\la)=0$.

Finally, since $\Omega(\ol)^*=-\Omega(\la)$ the first row of $\Omega$ is also zero.
\end{proof}

At present we cannot show that the remaining entries of $\Omega$ will also vanish under all circumstances even though we strongly suspect that this is the case.
However, we can show that it is true when $P_\pm(\la)=\id$ for all $\la\in\varrho(T)\setminus\tilde\Lambda$ (this includes all regular problems), when $n=1$, and when $n=2$ with $J$, $q$ and $w$ real.

\subsubsection{The case when \texorpdfstring{$P_\pm(\la)=\id$}{P+/-=1}}
\begin{lem}\label{L.4.4}
If $\la\in\varrho(T)\setminus\tilde\Lambda$ and $P_\pm(\la)=\id$, then $\Omega(\la)=0$ and hence $\bb M(\la)=\bb M(\ol)^*$.
\end{lem}

\begin{proof}
Only $\Omega_{4,4}$ needs to be considered in this case.
Recalling \eqref{wronski} we obtain
$$\Omega_{4,4}(\la)=(g^*Jg)^-(b)-(g^*Jg)^+(a)$$
which vanishes by Theorem \ref{T:3.5} since we have a self-adjoint restriction of $T_\mx$.
\end{proof}

\subsubsection{The case when $n=1$}
\begin{lem}
If $\la\in\varrho(T)\setminus\tilde\Lambda$ and $n=1$, then $\Omega(\la)=0$ and hence $\bb M(\la)=\bb M(\ol)^*$.
\end{lem}

\begin{proof}
First assume, by way of contradiction, that $\Omega_{2,2}(\la)\neq0$, i.e., $P_-(\la)=P_-(\ol)=0$.
Consider the problem where $\tilde q=q\chi_{(a,x_1)}$ and $\tilde w=w\chi_{(a,x_1)}$.
We now have that $b$ is a regular endpoint but the corresponding deficiency indices $\tilde n_\pm$ are still zero.
In other words we have a self-adjoint situation with $\tilde N=0$, $\tilde Q_-=1$, $\tilde Q_+=0$, and $\tilde{\bb P}=1$.
In this situation we find $\tilde{\bb M}(\la)=J^{-1}/2$ in either half plane.
But $J^{-1}$ is a non-zero purely imaginary number and, according to Lemma \ref{L.4.1}, $\tilde{\bb M}(\la)=\tilde{\bb M}(\ol)^*$, a contradiction.
Thus $\Omega_{2,2}(\la)=0$ and one shows similarly that $\Omega_{3,3}(\la)=0$.
We also have, trivially, that $\Omega_{3,2}(\la)$ and $\Omega_{2,3}(\la)$ are $0$.

We now consider the fourth row of $\Omega$.
The entry $\Omega_{4,2}(\la)$ contains the factor $P_-(\la)J^{-1}{(1-P_-(\la))}$ which is $0$ since $J^{-1}$ is a scalar.
It follows similarly $\Omega_{4,3}(\la)$ vanishes.
Finally, since $g(\cdot,\la)$ vanishes near $a$ or $b$ when $P_-(\la)=0$ or $P_+(\la)=0$, respectively, we may remove the factors $P_\pm$ in the expression
\begin{multline*}
\Omega_{4,4}(\la)
=-\big(g^*JU_0(\cdot,\lambda)P_-(\lambda)J^{-1}P_-(\ol)U_0(\cdot,\ol)^*Jg^*\big)^+(a)\\
 +\big(g^*JU_N(\cdot,\lambda)P_+(\lambda)J^{-1}P_+(\ol)U_N(\cdot,\ol)^*Jg^*\big)^-(b).
\end{multline*}
Thus, using \eqref{wronski} and \eqref{200820.1}, $\Omega_{4,4}(\la)=(g^*Jg)^-(b)-(g^*Jg)^+(a)=0$.
\end{proof}

\subsubsection{The case when $n=2$ and the coefficients are real}
The condition that $J$ is real, skew-adjoint, and invertible implies that $J=\beta \sm{0&-1\\ 1&0}$ for some $\beta\in \bb R\setminus\{0\}$.
Without loss of generality we shall henceforth assume that $\beta=1$.%
\footnote{If $\pm\beta>0$ we could employ the Liouville transform $v(x)=\pm\beta u(\pm x)$ to arrive at an equation of the same character but with $J=\sm{0&-1\\ 1&0}$.}
Next note that $u$ solves $Ju'+qu=\la wu$ if and only if $\ov u$ solves $J\ov u'+q\ov u=\ol w\ov u$.
This implies $n_+=n_-$ and $P_\pm(\ol)=\ov{P_\pm(\la)}$.

\begin{lem}
If $\la\in\varrho(T)\setminus\tilde\Lambda$, $n=2$, and $J$, $q$, and $w$ are real, then $\Omega(\la)=0$ and hence $\bb M(\la)=\bb M(\ol)^*$.
\end{lem}

\begin{proof}
As in the case $n=1$ we will make use of the auxiliary problem where $\tilde q=q\chi_{(a,x_1)}$ and $\tilde w=w\chi_{(a,x_1)}$.
Again we have that $b$ is a regular endpoint and that $\tilde N=0$.
We will also denote other quantities associated with the auxiliary problem by adding the $\tilde{}$ symbol.

First we show that it is impossible to have $P_\pm(\la)=0$.
Assuming, by way of contradiction, that $P_-(\la)=0$ and hence $P_-(\ol)=0$, we have $n_\pm=0$ and obtain $\tilde{\bb F}(\la)=(\id,0,0)^\top$ and $\tilde{\bb H}(\la)=\tilde{\bb F}(\la)/2$ for the auxiliary problem.
Consequently $\tilde{\bb M}(\la)=J^{-1}/2$ in either half-plane, which is, as before, absurd.
A similar argument shows that $P_+(\la)$ cannot be $0$.

When $P_-(\la)=\id$ then $\ms Q_-(\la)=0$ which implies that the second column and the second row of $\Omega(\la)$ vanish.

Next assume that the rank of $P_-(\la)$ is $1$.
Any orthogonal projection $P$ of rank $1$ in $\bb C^2$ satisfies $PJ^{-1}\overline P=0$ as a direct computation shows.
Applying this to $P=\id-P_-(\la)$ shows that $\Omega_{2,2}(\la)=0$.
That $\Omega_{3,2}(\la)=0$ is again trivial.
Now consider $\ms{A}_-(\lambda)=-(g^*J U_0(\cdot,\la)P_-(\la))^+(a)$ and let $\eta$ be an arbitrary element of $\bb C^2$.
If $u=U_0(\cdot,\la)P_-(\la)\eta$ and $\int_{(a,x_1)}u^*wu>0$, then our auxiliary problem is definite and we obtain $(g^*Ju)^+(a)=0$ and hence $\ms{A}_-(\lambda)=0$ from Lemma 7.6. in \cite{MR4047968}.
If $\int_{(a,x_1)}u^*wu=0$, i.e., if $u\in\tilde{\mc L}_0$, then $u$ satisfies $Ju'+\tilde qu=\mu\tilde w u=0$ for any $\mu\in\bb C$.
In this case there is an $\alpha\in\bb R^2$ such that $\tilde{\mc L}_0$ is spanned by $v=U_0(\cdot,0)\alpha$.
Both $u$ and $g$ are then multiples of $v$ and $(g^*Ju)^+=c\alpha^* U_0^+(\cdot,0)^*J U_0^+(\cdot,0)\alpha$ for a suitable $c\in\bb C$.
Because of equation \eqref{wronski} and since $\alpha^*J\alpha=0$, we obtain $\ms{A}_-(\lambda)=0$.
It follows that the second column and the second row of $\Omega(\la)$ vanish also when $P_-(\la)$ has rank $1$, i.e., in any case.

We may prove similarly that the third column and the third row of $\Omega(\la)$ vanish whatever $P_+(\la)$ may be and it remains only to consider
$$\Omega_{4,4}(\la)=\ms{A}_-(\lambda)\mathcal{J}^{-1}\ms{A}_-(\ol)^*-\ms A_+(\lambda)\mc J^{-1}\ms A_+(\bar{\lambda})^*.$$
If $P_-(\la)$ has rank $1$ we have $\ms{A}_-(\lambda)\mathcal{J}^{-1}\ms{A}_-(\ol)^*=0=(g^*J g)^+(a)$.
If $P_-(\la)=\id$ we get $\ms{A}_-(\lambda)\mathcal{J}^{-1}\ms{A}_-(\ol)^*=(g^*J g)^+(a)$ on account of equation \eqref{wronski}.
With similar considerations for $P_+(\la)$ we have therefore $\Omega_{4,4}(\la)=(g^*J g)^+(a)-(g^*J g)^+(b)$
which is $0$ by \eqref{200820.1}.
\end{proof}

Later we will also need the following lemma.
\begin{lem}\label{L:4.7}
If $\mathbf{B}=\ran\bb P$ or if $\Omega(\la)=0$, then $\ran\bb H(\la)\mc J^{-1}\bb P\subset\ran\bb F(\la)$ and $\ran\bb F(\la)^\dag\bb H(\la)\mc J^{-1}\bb P\subset\ran\bb P$.
\end{lem}

\begin{proof}
If $\mathbf{B}=\ran\bb P$ and $z$ is any element of $\bb C^{n(N+1)}$, then there is a compactly supported $f\in L^2(w)$ such that $\bb Pz =(\ms Ff)(\la)$.
Therefore \eqref{201106.1} and \eqref{221107.1} establish the existence of a vector $\tilde z=\frac12(\tilde u+\tilde v)\in\ran\bb P$ such that $\bb F(\la)\tilde z
=\bb H(\la)\mc J^{-1}\bb Pz$ settling our first claim in this case.
The second follows after applying $\bb F(\la)^\dag$ since $\bb F(\la)^\dag\bb F(\la)=\id$.

If $\Omega(\la)=0$ we have that  $\mathbb{H}(\lambda)\mathcal{J}^{-1}\mathbb{P}\mathbb{F}(\ol)^* =-\mathbb{F}(\lambda)\mathbb{P}\mathcal{J}^{-1}\mathbb{H}(\ol)^*$.
Here we use $\bb F(\ol)^*\bb F(\ol)^{\dag *}=\id$ to get $\mathbb{H}(\lambda)\mathcal{J}^{-1}\mathbb{P}=-\mathbb{F}(\lambda)\mathbb{P}\mathcal{J}^{-1}\mathbb{H}(\ol)^*\mathbb{P}\mathbb{F}(\ol)^{\dag*}$.
\end{proof}

\begin{remark}
If $\mathbf{B}=\ran\bb P$ we have now that $\ran\bb F(\la)^\dag\Omega(\la)\subset\ran\bb P$.
Since $\bb F(\la)\bb F(\la)^\dag$ is the orthogonal projection onto the range of $\bb F(\la)$, this shows that
$\bb F(\la)(\bb M(\la)-\bb M(\ol)^*)\bb F(\ol)^*=\Omega(\la)=0$.
In other words, the requirement $\Omega(\la)=0$ is satisfied when $\mathbf{B}=\ran\bb P$.
\end{remark}

\subsection{The imaginary part of \texorpdfstring{$\bb M$}{M}}
In this section we assume that $\la$ is in the upper half-plane but not in $\tilde\Lambda$ and that $\Omega(\la)=0$ (which holds when $\mathbf{B}=\ran\bb P$).
Let $\mc S(x)$ be the diagonal block matrix whose entries are, in this order, the $n\times n$ blocks $\sgn(x-\xi_j)\id$, $j=0,...,N$.
Recall from the beginning of Section \ref{S:sol} that $U_j(\xi_j,\la)=\id$.

Define $\theta$ by
$$\theta(x,\la)=\big(\bb M(\la)+\frac12 \mc S(x)\mc J^{-1}\bb P\big) z$$
with $z\in\bb C^{n(N+1)}$.
Note that $\ms U(\cdot,\la)\theta(\cdot,\la)$ is a solution of $Ju'+(q-\la w)u=0$ in each of the intervals $(x_k,\xi_k)$ and $(\xi_k,x_{k+1})$, $k=0,...,N$.
We want to show that it actually satisfies the differential equation on $(\xi_{k-1},\xi_k)$, $k=1,...,N$, as well as the boundary conditions including the requirement that $\ms U(\cdot,\la)\theta(\cdot,\la)$ is in $\mc L^2(w)$.
For the first claim we need to show that
\begin{equation}\label{221112.1}
\bb F_k(\la)\theta(x_k,\la)=0
\end{equation}
where $\bb F_k$ denotes the $k$-th row of $n\times n$ blocks of $\bb F$ and hence of $\bb B$.
Using $\bb F_k\bb P=\bb F_k$, see Lemma~\ref{L:3.1}, and $\bb F^\dag\bb F=\id$, the left-hand side of \eqref{221112.1} becomes
$$\bb F_k(\la)\bb F(\la)^\dag\big(\bb H(\la)+\frac12\bb F(\la)S(x_k)\big)\mc J^{-1}\bb P z.$$
Since $\bb F(\la)\bb F(\la)^\dag$ is the orthogonal projection onto $\ran\bb F(\la)$ and since, by Lemma~\ref{L:4.7}, $\bb H(\la)\mc J^{-1}\bb Pz$ is in $\ran\bb F(\la)$ we get next
$$\frac12\big(\tilde{\bb B}_k(\la)+\bb B_k(\la)S(x_k)\big)\mc J^{-1}\bb P z.$$
This does indeed vanish due to the special structure of $\bb B$ and $\tilde{\bb B}$ and since $\mc S(x_k)$ is a diagonal block matrix whose first $k$ blocks are $\id_n$ while the remaining $N+1-k$ blocks are $-\id_n$.

To show that $\ms U(\cdot,\la)\theta(\cdot,\la)$ is in $L^2(w)$ we are following a very similar strategy.
On the left-hand side of \eqref{221112.1} we have to choose $k=N+1$ so that $\bb F_k=\ms Q_-$ or $k=N+2$ so that $\bb F_k=\ms Q_+$.
We also have to choose $x$ in either $(a,\xi_0)$ or else in $(\xi_N,b)$ so that $S(x)$ is either $-\id_{n(N+1)}$ or else $\id_{n(N+1)}$.
Since $\ms Q_\pm\bb P=\ms Q_\pm$ we get $\bb F_k(\la)\bb M(\la)=\bb H_{k}(\la)\mc J^{-1}\bb P=\pm\frac12 \bb F_k(\la)\mc J^{-1}\bb P$.
This proves \eqref{221112.1} for $k=N+1$ and $k=N+2$.

Finally, the boundary condition translates to $\ms A_+(\la)\theta_+(\la)+\ms A_-(\la)\theta_-(\la)=0$ where $\theta_\pm(\la)=\theta(x,\la)$ with $x>\xi_N$ for the upper sign and $x<\xi_0$ for the lower sign.
Note that, imitating previous arguments,
$$(\ms A_+(\la)+\ms A_-(\la))\theta_+(\la)=(\bb H_{N+3}(\la)+\frac12 \bb F_{N+3}(\la))\mc J^{-1}\bb Pz=\ms A_-(\la)\mc J^{-1}\bb Pz$$
gives $\ms A_+(\la)\theta_+(\la)=\ms A_-(\la)(\frac12\mc J^{-1}\bb Pz-\bb M(\la)z)=-\ms A_-(\la)\theta_-(\la)$, our desired result.

We now abbreviate $\ms U\theta$ by $s$.
Fix $\la,\mu\in\varrho(T)\setminus\tilde{\Lambda}$ and define $h=s(\cdot,\la)-s(\cdot,\mu)$.
Thus $h$ satisfies the differential equation
\[Jh'+(q-\mu w)h=(\la-\mu)ws(\cdot,\la)\]
in the intervals $(a,\xi_0)$, $(\xi_{k-1},\xi_k)$, $k=1,...,N$, and in $(\xi_N,b)$.
Since the jump of $s(\cdot,\la)$ at $\xi_k$ is equal to the $k$-th row of blocks $\mc J^{-1}\bb Pz$ and hence independent of $\la$ it follows that $h$ is continuous at those points which entails that it satisfies the above equation in all of $(a,b)$.
Since $h$ is also in $\mc L^2(w)$ and satisfies the boundary conditions (if any) it is thus the case that $h$ is in the class of $(\la-\mu)R_\mu s(\cdot,\la)$.
Evaluating $h$ at $\xi_{k-1}$ for $ k=1,...,N+1 $ gives
\begin{equation}\label{220414}
(h(\xi_0),...,h(\xi_n))^\diamond=(\mathbb{M}(\la)-\mathbb{M}(\mu))z
\end{equation}
showing that $(h(\xi_0),...,h(\xi_n))^\diamond\in N_0^\perp$ so that
$h=(\la-\mu)\mathcal{R}_\mu s(\cdot,\la)$.

Multiply \eqref{220414} by $z^*$ on the left to get
\begin{equation}\label{221108.1}
\frac{z^*\mathbb{M}(\la)z-z^*\mathbb{M}(\mu)z}{\la-\mu}=\sum_{k=1}^{N+1}z_k^* (\mathcal{R}_\mu s(\cdot,\la))(\xi_{k-1}).
\end{equation}
Using \eqref{elf} and $\ms U(\xi_{k-1},\mu)=e_k^*$, where $e_k^*$ is a row of $n\times n$ blocks all zero except for the $k$-th which is $\id_n$, we get
$$(\mc R_\mu s(\cdot,\la))(\xi_{k-1})=e_k^*\int \big(\bb M(\mu)+\frac12\mc J^{-1}\sgn(\xi_{k-1}-\cdot)\big)\ms U(\cdot,\overline\mu)^*w s(\cdot,\la).$$
Now let $\mu=\ol$ and recall that $\bb M(\ol)=\bb M(\la)^*$.
Then
$$z^*\frac{\mathbb{M}(\la)-\mathbb{M}(\la)^*}{\la-\ol}z
 =\int \sum_{k=1}^{N+1}\big[\ms U(\cdot,\la)\big(\bb M(\la)-\frac12\mc J^{-1}\sgn(\xi_{k-1}-\cdot)\big)e_kz_k\Big]^*w s(\cdot,\la).$$

Using the identities $\sum_{k=1}^{N+1}e_kz_k=z$ and $\sum_{k=1}^{N+1}\sgn(\xi_{k-1}-x)e_kz_k=-\mc S(x)z$ and the fact that the matrices $\mc J^{-1}$ and $\mc S(x)$ commute, we get now
$$z^*\frac{\mathbb{M}(\la)-\mathbb{M}(\ol)}{\la-\ol}z=\int s(\cdot,\la)^*ws(\cdot,\la)$$
which is non-negative.

We have proved the following lemma.
\begin{lem}
If $\Im\la>0$, $\la\notin\tilde\Lambda$ and either $\mathbf{B}=\ran\bb P$ or $\Omega(\la)=0$, then $\Im M(\la)\geq0$.
\end{lem}

\subsection{Analyticity}
Suppose $\la$ is a non-real complex number for which there is a neighborhood which does not intersect $\tilde\Lambda$.
The resolvent identity for $R_\mu$ and the boundedness of the operator selecting a representative of $([u],[f])\in T_\mx$ (see Lemma 5.2 in \cite{MR4431055}) show that the map $\mu \mapsto \mc R_\mu f$ is continuous at $\la$.
Using now equation \eqref{221108.1} shows that $z^*(\bb M(\mu)-\bb M(\la))z/(\mu-\la)$ has a limit as $\mu$ tends to $\la$.
Thus $z^*\bb Mz$ is analytic near $\lambda$.

\subsection{\texorpdfstring{$\bb M$}{M} is Nevanlinna}
In the course of our investigations we have added several hypotheses to the basic Hypothesis \ref{hyp:m}.
We will add one more and collect them in the following Hypothesis \ref{hyp:fe}.

\begin{hyp}\label{hyp:fe}
$(a,b)$ is a real interval.
$J$ is a constant, invertible, skew-hermitian $n\times n$-matrix.
$q$ and $w$ are $n\times n$-matrices whose entries are distributions of order $0$ on $(a,b)$; $q$ is hermitian and $w$ is non-negative.
$\Lambda_x\cap\bb R$ is empty unless $x\in\{x_1,...,x_N\}\subset (a,b)$ and $\tilde\Lambda=\bigcup_{x\not\in\{x_1,...,x_N\}}\Lambda_{x}$ is a closed set of isolated points.
Finally, we require $\Omega=0$ (which is satisfied when $\mathbf{B}=\ran\bb P$).
\end{hyp}

Under this hypothesis we can prove the following theorem.
\begin{thm}
Assume the validity of Hypothesis \ref{hyp:fe}.
Then the function $\bb M$ may be extended to all of $\bb C\setminus\bb R$ as a matrix-valued Nevanlinna function.
\end{thm}

\begin{proof}
Let $z\in\mathbb{C}^{n(N+1)}$ and $m=z^*\mathbb{M}z$.
The singularities of $m$ are the points in $\tilde\Lambda$, a closed set of isolated points by hypothesis.
Suppose now that $\mu$ is one of these points and that $B$ is a ball centered at $\mu$ not intersecting $\bb R$ or $\tilde\Lambda\setminus\{\mu\}$.
Note that $\Im m(\lambda)/\Im\lambda\geq0$ in $B\setminus\{\mu\}$ and hence $\mu$ is a removable singularity of $m$.
Since this is so for any $z\in\bb C^{n(N+1)}$ we may extend $\bb M$ to all $\bb C\setminus\bb R$ as an analytic function.
The properties of symmetry and of the sign of $\Im \bb M$ are also retained.
\end{proof}

\section{The Fourier expansion}\label{S:FE}
We begin with a few words about the spectral theorem for self-adjoint relations $T$.
The closure $\mc H_0$ of the domain of $T$ is the orthogonal complement of $\mc H_\infty=\{f\in L^2(w):(0,f)\in T\}$ in $L^2(w)$.
It follows that $T=T_0\oplus (\{0\}\times \mc H_\infty)$ where $T_0=T\cap(\mc H_0\times\mc H_0)$ is a self-adjoint operator densely defined in $\mc H_0$.
The spectral theorem for self-adjoint operators guarantees the existence of a resolution of the identity $\pim$ such that
$$\<f,T_0g\>=\int t\;d\<f,\pim((-\infty,t))g\>.$$
One may now extend the domain of definition of the spectral projections $\pim(B)$ from $\mc H_0$ to $\mc H$ by setting $\pim(B)f=0$ whenever $f\in \mc H_\infty$.
Thus $\pim(\bb R)$ becomes the orthogonal projection from $\mc H$ onto $\mc H_0$.
For more information we refer the reader to Appendix B of \cite{MR4047968} or to Section 2.2 of \cite{MR4199125}.

In the following we require Hypothesis \ref{hyp:fe} to hold.
Then, as we just proved $\bb M$ is a Nevanlinna function and therefore defines a spectral measure $\tau$, see equation \eqref{mnl}.
For compactly supported functions $f$ of $L^2(w)$ and if $\la\in\bb C\setminus\tilde\Lambda$ we defined earlier the Fourier transform $\ms F$ by
$$(\ms F f)(\la)=\int \ms U(\cdot,\ol)^*wf.$$
Note that $\bb C\setminus\tilde\Lambda$ is an open set containing the real line, so $\ms Ff$ is defined everywhere in a neighborhood of $\bb R$.
In fact, we are mostly concerned with the restriction of $\ms Ff$ to $\bb R$ but will not use different notation.

The remainder of this section is devoted to the proof of the following theorem, the main theorem of this paper.
The outline of the proof follows, with one exception, the proof of Theorem 15.5 in Bennewitz \cite{Bennewitz-sths}.
The spirit of that proof was also used in \cite{MR4199125} and in \cite{MR4047968} which may be consulted for some of the details we skip here in the interest of brevity.
The exception concerns the proof of the fact that $\ker \ms F^*$ is trivial, see Lemma ~\ref{kerg} below where $\ms F^*$ is called $\ms G$.

\begin{thm} \label{t:main}
Suppose $T$ is a self-adjoint restriction of a relation $T_\mx$ whose coefficients $q$ and $w$ satisfy Hypothesis \ref{hyp:fe}.
Let $\tau$ be the measure generated by the associated $\bb M$-function.
Then the following statements hold.
\begin{enumerate}
\item There is a continuous map $\ms F:L^2(w)\to L^2(\tau)$ which assigns to a compactly supported element $f\in L^2(w)$ the function defined by $(\ms Ff)(t)=\int \ms U(\cdot,t)^*wf$.
The kernel of $\ms F$ is the space $\mc H_\infty=\{f\in L^2(w):(0,f)\in T\}$.
\item  There is a continuous map $\ms G:L^2(\tau)\to L^2(w)$ which assigns to a compactly supported element $\hat f\in L^2(\tau)$ the function defined by $(\ms G\hat f)(x)=\int \ms U(x,\cdot)\tau \hat f$.
The range of $\ms G$ is the space $\mc H_0=\ov{\dom T}$.
\item $\ms F\circ\ms G$ is the identity operator on $L^2(\tau)$, $\ms G\circ\ms F$ is the orthogonal projection from $L^2(w)$ onto $\mc H_0$, and the restriction of $\ms F$ to $\mc H_0$ is unitary.
\item If $(u,f)\in T$ then $(\ms Ff)(t)=t (\ms Fu)(t)$.
Conversely, if $t\mapsto \hat u(t)$ and $t\mapsto \hat f(t)=t\hat u(t)$ are both in $L^2(\tau)$, then $(\ms G\hat u,\ms G\hat f)\in T$.
\end{enumerate}
\end{thm}

\begin{remark}
When $q$ and $w$ are finite measures on $(a,b)$, i.e., when the endpoints $a$ and $b$ are regular, the conditions that $\tilde\Lambda$ is a closed set of isolated points and that $\Omega=0$ in Hypothesis \ref{hyp:fe} are automatically satisfied.
\end{remark}

Since the functions $\ms U(x,\cdot)$ are analytic on $\bb C\setminus\tilde\Lambda$ it follows immediately that
$$\<g,R_\la f\>-(\ms Fg)(\ol)^*\bb M(\la)(\ms Ff)(\la)$$
extends to an analytic function on $\bb C\setminus\tilde\Lambda$ so that we obtain the identity
\begin{equation}\label{220810.1}
\oint_\Gamma \<f,R_\la f\>\,d\la=\oint_\Gamma(\ms Ff)(\ol)^* \bb M(\la)(\ms Ff)(\la)\,d\la
\end{equation}
when $\Gamma$ is the contour described by the rectangle with vertices $c\pm i\varepsilon$ and $d\pm i\varepsilon$ for $c<d$ and sufficiently small but positive $\varepsilon$, provided the integrals exist.

Let $\pim$ be the resolution of the identity for our self-adjoint relation $T$ and define
$$\Pim_{f,g}(t)=\<f,\pim((-\infty,t))g\>_w.$$
With the aid of the spectral theorem, Fubini's theorem, and Cauchy's integral formula one shows now that
\begin{equation}\label{180125.2}
\oint_\Gamma \<f,R_\la f\>d\la=-2\pi i \int_{[c,d)} d\Pim_{f,f}.
\end{equation}
when $c$ and $d$ are points of differentiability of $\Pim_{f,f}$.
Similarly, using the Nevanlinna representation of $\bb M$, Fubini's theorem, and Cauchy's integral formula gives
\begin{equation}\label{180125.1}
\oint_\Gamma(\ms Ff)(\ol)^*\bb M(\la)(\ms Ff)(\la)d\la=-2\pi i \int_{[c,d)} (\ms Ff)^*\tau(\ms Ff)
\end{equation}
when $c$ and $d$ are points of differentiability of an antiderivative $\mathfrak{T}$ of $\tau$.
Thus we obtain
\begin{equation}\label{180212.1}
\int_{[c,d)} (\ms Ff)^*\tau(\ms Ff)=\int_{[c,d)} d\Pim_{f,f}.
\end{equation}
Since $\Pim_{f,f}$ is left-continuous and $\mathfrak{T}$ may be chosen to be, equation \eqref{180212.1} actually holds for all $c$, $d$ with $c<d$.

With these preparations we may prove statement (1) of Theorem \ref{t:main}.
\begin{lem}\label{L:6.6}
When $f\in L^2(w)$ is compactly supported, then $\ms Ff$ is in $L^2(\tau)$.
The map $\ms F$ extends, by continuity, to all of $L^2(w)$.
Moreover,
\begin{equation}\label{221113.1}
\Pim_{f,g}(t)=\<f,\pim((-\infty,t))g\>_w=\int_{(-\infty,t)} (\ms F f)^* \tau (\ms Fg)
\end{equation}
whenever $f,g\in L^2(w)$.
In particular, $\<f,\pi(\bb R)g\>_w=\<\ms F f,\ms F g\>_\tau$ and $\ker\ms F=\ms H_\infty$.
\end{lem}

\begin{proof}
Suppose $f\in L^2(w)$ is compactly supported.
Choosing the interval $[c,d)$ sufficiently large in equation \eqref{180212.1} proves that $\ms Ff\in L^2(\tau)$.
If $f$ is an arbitrary element of $L^2(w)$ and $n\mapsto[a_n,b_n]$ a sequence of intervals in $(a,b)$ converging to $(a,b)$, set $f_n=f\chi_{[a_n,b_n]}$.
Then
$$\|\ms F f_n-\ms Ff_m\|_\tau=\|\pi(\bb R)(f_n-f_m)\|_w\leq \|f_n-f_m\|_w$$
showing that $n\mapsto \ms Ff_n$ is a Cauchy sequence in $L^2(\tau)$ and thus convergent.
By interweaving sequences it follows that the limit of this Cauchy sequence does not depend on how $f$ is approximated.
We denote the limit by $\ms F f$ thereby extending our definition of the Fourier transform to all of $L^2(w)$.
Equation \eqref{221113.1} holds when $g=f$ and otherwise by polarization.
\end{proof}

Now we define a transform $\ms G:L^2(\tau)\to L^2(w)$.
We begin by setting
$$(\ms G\hat f)(x)=\int_{\bb R} \ms U(x,\cdot)\tau\hat f$$
whenever $\hat f$ is compactly supported.
Note that $\ms G\hat f$ is locally of bounded variation.
Then we have the following result which proves statement (2) of Theorem \ref{t:main}.

\begin{lem}\label{L:5.8}
When $\hat f\in L^2(\tau)$ is compactly supported, then $\ms G\hat f$ is in $L^2(w)$.
The map $\ms G$ extends, by continuity, to all of $L^2(\tau)$.
We have that
\begin{equation}\label{180110.1}
\<g,\ms G\hat f\>_w=\<\ms Fg,\hat f\>_\tau
\end{equation}
for all $g\in L^2(w)$ and all $\hat f\in L^2(\tau)$.
Moreover, $\ker\ms G=(\ran\ms F)^\perp$, $\ran\ms G=H_0$, and $\ms G\circ\ms F=\pi(\bb R)$.
\end{lem}

\begin{proof}
Suppose that $\hat f\in L^2(\tau)$ is compactly supported, denote $\ms G\hat f$ by $f$, and let $f_n=f\chi_{[a_n,b_n]}$.
Upon changing the order of integration we get
$$\|f_n\|_w^2=\<f_n,f\>_w=\<\ms Ff_n,\hat f\>_\tau\leq \|\ms Ff_n\|_\tau \|\hat f\|_\tau.$$
Lemma \ref{L:6.6} implies $\|\ms Ff_n\|_\tau=\|\pi(\bb R)f_n\|_w\leq \|f_n\|_w$ and hence $\|f_n\|_w\leq \|\hat f\|_\tau$.
This is the case for every interval $[a_n,b_n]\subset (a,b)$ so it follows that $\ms G\hat f\in L^2(w)$.

As before we extend the domain of definition of $\ms G$ from the compactly supported functions in $L^2(\tau)$ to all of $L^2(\tau)$.
Specifically, for a general element $\hat f$ in $L^2(\tau)$ set $\hat f_n=\hat f\chi_{[-n,n]}$.
Then, according to what we just proved, $\|\ms G\hat f_n-\ms G\hat f_m\|_w\leq \|\hat f_n-\hat f_m\|_\tau$ which implies that $n\mapsto \ms G\hat f_n$ is a Cauchy sequence in $L^2(w)$ and thus convergent.
We denote the limit by $\ms G\hat f$.

Now suppose $g\in L^2(w)$, $\hat f\in L^2(\tau)$, $[a_k,b_k]\subset(a,b)$ and $[-n,n]\subset\bb R$.
Upon changing the order of integration we get (as before)
$$\<g\chi_{[a_k,b_k]},\ms G (\hat f\chi_{[-n,n]})\>_w=\<\ms F(g\chi_{[a_k,b_k]}),\hat f\chi_{[-n,n]}\>_\tau.$$
Now let $[a_k,b_k]\times[-n,n]$ approach $(a,b)\times\bb R$ to obtain equation \eqref{180110.1}.

Equation \eqref{180110.1} implies $\ker\ms G=(\ran\ms F)^\perp$ and $\ran\ms G\subset \mc H_\infty^\perp=\mc H_0$.
Choosing $\hat f=\ms F f$ in \eqref{180110.1} implies, using Lemma \ref{L:6.6}, that $\<g,(\ms G\circ\ms F)f\>=\<g,\pi(\bb R)f\>$.
Since this is so for all $g\in L^2(w)$ we get $\ms G\circ\ms F=\pi(\bb R)$ and, in particular $\mc H_0\subset\ran\ms G$.
\end{proof}

Our next goal is to show that $\ker\ms G$ is trivial so that $\ran\ms F$ is dense in $L^2(\tau)$.
This will show that statement (3) of Theorem \ref{t:main} holds as can be seen as follows.
Using Lemma \ref{L:6.6}, Lemma \ref{L:5.8}, and the self-adjointness of $\pi(\bb R)$ we obtain
$$\<\ms Fg,(\ms F\circ\ms G-\id)\hat f\>=\<\ms Fg,\ms F(\ms G\hat f)\>-\<\ms Fg,\hat f\>=0.$$
This implies $\ms F\circ\ms G=\id$ and that $\ms G$ is an isometry.
We have already that $\ms G\circ\ms F=\pi(\bb R)$ and that $\ms F|_{\mc H_0}$ is an isometry.
However, the proof of $\ker\ms G=\{0\}$, see Lemma \ref{kerg}, requires more preparation.

\begin{lem}\label{L:6.7}
If $\Im(\la)\neq0$, then $(\ms F(R_\la g))(t)=(\ms Fg)(t)/(t-\la)$.
\end{lem}

\begin{proof}
First note that $t\mapsto \hat g(t)/(t-\la)$ is in $L^2(\tau)$ if $\hat g$ is.
The spectral theorem and Lemma \ref{L:6.6} give
$$\<f,R_\la g\>_w=\int \frac1{t-\la}\ d\Pim_{f,g}(t)=\int \frac{(\ms F f)(t)^* \tau(t) (\ms Fg)(t)}{t-\la}=\<\ms Ff,(\ms Fg)/(\cdot-\la)\>_\tau.$$
In particular, $\|R_\la g\|_w^2=\<\ms F(R_\la g),(\ms Fg)/(\cdot-\la)\>_\tau$.
On the other hand
$$\|R_\la g\|_w^2=\<g,R_\ol R_\la g\>_w=\frac1{\la-\ol}\<g,(R_\la- R_\ol)g\>_w=\|\ms Fg/(\cdot-\la)\|^2_\tau.$$
Lemma \ref{L:6.6} also implies that $\|R_\la g\|_w^2=\|\ms F(R_\la g)\|_\tau^2$.
Thus the four terms appearing in the expansion of $\|\ms F(R_\la g)-\ms Fg/(\cdot-\la)\|^2$ cancel each other leaving $0$.
\end{proof}

Let $\bb T=(\tau/\tr\tau)$ be the Radon-Nikodym derivative of $\tau$ with respect to $\tr\tau$.
Note that $\bb T\in L^\infty(\tr\tau)$.
By the Lebesgue-Radon-Nikodym theorem we have $\tr\tau=h \bm+\sigma$ where $\bm$ denotes Lebesgue measure, $h$ is a non-negative function, $\sigma$ is a non-negative measure, and $\sigma$ and $\bm$ are mutually singular.
Define
$$\omega(s,\varepsilon)=\int_{\bb R} \frac{\varepsilon}{(s-t)^2+\varepsilon^2} \tr\tau(t).$$
By Fatou's theorem (see, e.g., Theorem 5.5 in Rosenblum and Rovnyak \cite{MR1307384}) $\omega(s,\varepsilon)$ converges to $\pi h(s)$ for $\bm$-almost every $s\in\bb R$ as $\varepsilon\downarrow0$.
In fact $h(s)>0$ $\bm$-almost everywhere.
The measure $\sigma$ is concentrated on the set
$$S=\{s\in\bb R: \lim_{r\downarrow0} \tr\tau((s-r,s+r))/(2r)=\infty\}$$
and, consequently, when $s\in S$ then $\omega(s,\varepsilon)$ tends to $\infty$ as $\varepsilon\downarrow0$.
It follows that $1/\omega(s,\varepsilon)$ is bounded above when $s$ is in a set of full $\tr\tau$-measure.

\begin{lem}\label{bt}
$\bb B\bb T=0$ and $(1-\bb P)\bb T=0$ almost everywhere with respect to $\tr\tau$.
\end{lem}

\begin{proof}
Since $\bb F(\la)\bb F(\la)^\dag$ is the orthogonal projection onto the range of $\bb F(\la)$ Lemma \ref{L:4.7} shows that
$$\bb H(\la)\mc J^{-1}\bb P = \bb F(\la)\bb F(\la)^\dag\bb H(\la)\mc J^{-1}\bb P.$$
Using that $\bb B=\bb B\bb P$ (see Lemma \ref{L:3.1}) the first $nN$ rows of this identity are
$$\bb B(\la)\bb M(\la)=\frac12\tilde{\bb B}(\la)\mc J^{-1}\bb P.$$
Subtract from this identity the one where $\la$ is replaced by $\ol$ to get
\begin{equation}\label{220624.1}
\bb B(\la)\bb M(\la)-\bb B(\ol)\bb M(\ol)=\frac12(\tilde{\bb B}(\la)-\tilde{\bb B}(\ol))\mc J^{-1}\bb P.
\end{equation}
We will now take the limit as $\varepsilon=\Im\la\downarrow0$.
First note that the right-hand side of \eqref{220624.1} will tend to $0$.
Next, the integral occurring in $\iu\varepsilon \bb M(s\pm\iu\varepsilon)$ is
$$\int_{\bb R}\frac{\iu\varepsilon (1+ts\pm\iu t\varepsilon)}{t-s\mp\iu\varepsilon} \frac{\tau(t)}{t^2+1}.$$
For $\varepsilon\in[0,1]$ and $t\in\bb R$ the first fraction may be bounded by $5(s^2+1)$.
Since the measure $\tau(t)/(t^2+1)$ is finite, the dominated convergence theorem shows that
$$\lim_{\varepsilon\downarrow0}\iu\varepsilon \bb M(s\pm\iu\varepsilon)=\mp\tau(\{s\}).$$
Also $\lim_{\varepsilon\downarrow0}(\bb B(s\pm\iu\varepsilon)-\bb B(s))/(\iu\varepsilon)=\pm\dot{\bb B}(s)$.
These facts may be combined to give
$$\lim_{\varepsilon\downarrow0}\bb B(s)(\bb M(\la)-\bb M(\ol))=0.$$
The Nevanlinna representation of $M$ gives
$$\bb M(s+\iu\varepsilon)-\bb M(s-\iu\varepsilon)=2\iu\varepsilon B+2\iu\int_{\bb R} \frac{\varepsilon}{(t-s)^2+\varepsilon^2}\bb T(t) \tr\tau(t).$$
Using the fact that the function $1/\omega(s,\varepsilon)$ is bounded above for $\tr\tau$-almost every $s$, and Theorem \ref{fatou} we get now
$$\bb B(s)\bb T(s)=\lim_{\varepsilon\downarrow0} \bb B(s)
 \int_{\bb R} \frac{\varepsilon}{(t-s)^2+\varepsilon^2}\bb T(t) \tr\tau(t)/\omega(s,\varepsilon)=0.$$

The proof of the identity $(1-\bb P)\bb T=0$ is very similar after remembering that $\bb P$ is the left-most factor of $\bb M$ and hence
$$0=(1-\bb P)(\bb M(\la)-\bb M(\ol))=2\iu(1-\bb P)\int_{\bb R} \frac{\varepsilon}{(t-s)^2+\varepsilon^2}\bb T(t) \tr\tau(t).$$
This completes the proof.
\end{proof}

\begin{lem}\label{kerg}
$\ker\mc G$ is trivial.
\end{lem}

\begin{proof}
Suppose $\hat u\in\ker\mc G$.
Then, using Lemmas \ref{L:5.8} and \ref{L:6.7}, we find
$$0=\<R_\ol g,\mc G\hat u\>_w=\<\ms F R_\ol g,\hat u\>_\tau=\int \frac{(\ms F g)(t)^*}{t-\la}\tau(t)\hat u(t)$$
for any $\la\in\bb C\setminus\bb R$ and any $g\in L^2(w)$.
The Stieltjes inversion formula shows that $(\ms Fg)^*\tau\hat u$ is the zero measure.
If $K$ is a compact subset of $\bb R$ and $\psi_K(t)=t\chi_K(t)$, then
$$\<g,\mc G\chi_K\hat u\>_w=\<\ms Fg,\chi_K\hat u\>_\tau=0 \disand \<g,\mc G\psi_K\hat u\>_w=\<\ms Fg,\psi_K\hat u\>_\tau=0$$
show that $\chi_K\hat u$ and $\psi_K\hat u$ are also in $\ker \mc G$.

Define $u(x)=\int_K \ms U(x,t)\bb T(t) \hat u(t)\tr\tau(t)$ and $f(x)=\int_K \ms U(x,t)\bb T(t) t\hat u(t)\tr\tau(t)$ which are representatives of $\mc G \chi_K\hat u=[0]$ and $\mc G \psi_K\hat u=[0]$, respectively.

We proved in Lemma \ref{bt} that $\bb B\bb T\hat u=0$ almost everywhere with respect to $\tr\tau$.
Therefore, by equation \eqref{eq:bb} for $f=0$, $\ms U(\cdot,t)\bb T(t)\hat u(t)$ is a solution of $Jv'+qv=twv$ on $(a,b)$ for almost every $t\in K$.
Using Fubini's theorem this implies that $Ju'+qu=wf$.
Since $u$ and $f$ are in $[0]$, it follows that $u\in\mc L_0$.
Therefore $\tilde u_0=(u(\xi_0),...,u(\xi_N))^\diamond=\int_K \bb T\hat u\tr\tau$ is in $\ker\bb P$.
By Lemma \ref{bt} we also have that $(\id-\bb P)\bb T=0$ and hence $\tilde u_0=0$.
This is only possible when $u$ is identically equal to $0$.

Choosing $K=[0,s]$ or $K=[s,0]$ shows that the cumulative distribution function of the measure $\ms U(x,\cdot)\bb T\hat u\tr\tau$ is zero, i.e., this measure is the zero measure.
In particular, if $x=\xi_k\in (x_{k},x_{k+1})$ this is the $(k+1)$-st block of $\bb T\hat u \tr\tau=\tau\hat u$.
But this means that $\hat u$ is zero almost everywhere with respect to $\tau$.
\end{proof}

Our last lemma provides the proof of statement (4) of Theorem \ref{t:main}.
\begin{lem}
If $(u,f)\in T$ then $(\ms Ff)(t)=t (\ms Fu)(t)$.
Conversely, if $t\mapsto \hat u(t)$ and $t\mapsto \hat f(t)=t\hat u(t)$ are both in $L^2(\tau)$, then $(\ms G\hat u,\ms G\hat f)\in T$.
\end{lem}

\begin{proof}
Suppose $(u,f)\in T$ and hence that $(f-\la u,u)\in R_\la$.
Then Lemma~\ref{L:6.7} gives $(\ms Fu)(t)=(\ms F(f-\la u))(t)/(t-\la)$ which simplifies to $t(\ms Fu)(t)=(\ms Ff)(t)$.

Now suppose that $\hat u,\hat f\in L^2(\tau)$ where $\hat f(t)=t\hat u(t)$.
Define $f=\ms G\hat f$ and $u=\ms G\hat u$ and pick a $\la$ in $\varrho(T)\setminus\tilde\Lambda$.
Then, using Lemma \ref{L:6.7},
$$\ms Fu=\hat u=\frac{\hat f-\la\hat u}{t-\la}=\frac{\ms F(f-\la u)}{t-\la}=\ms F(R_\la(f-\la u)).$$
Applying $\ms G$ gives $u=R_\la(f-\la u)$ which is equivalent to $(u,f)\in T$.
\end{proof}

\appendix
\section{An extension of Fatou's theorem}
In 1906 Fatou \cite{MR1555035} investigated the limiting behavior of holomorphic functions defined on the unit disk.
Analogous considerations for holomorphic functions defined on the upper half-plane lead to the following theorem, commonly called Fatou's theorem, see, e.g., Theorem 5.5 in Rosenblum and Rovnyak \cite{MR1307384}.

\begin{thm}
Let
$$V(s+\iu r)=\int_{\bb R}\frac{r\mu(t)}{(t-s)^2+r^2}$$
where $\mu$ is a non-negative Borel measure satisfying $\int_{\bb R} \mu(t)/(t^2+1)<\infty$, $r>0$, and $s\in\bb R$.
If the Lebesgue-Radon-Nikodym decomposition of $\mu$ is $h\bm+\sigma$ (where $\bm$ denotes Lebesgue measure), then
$$\lim_{r\downarrow0} V(s+\iu r)=\pi h(s)$$
almost everywhere with respect to Lebesgue measure.
\end{thm}

Applying this result to the measure $f\mu$ when $f\in L^\infty(\mu)$ gives
\begin{equation}\label{221115.1}
\lim_{r\downarrow0}\left[\int_{\bb R}\frac{rf(t)\mu(t)}{(s-t)^2+r^2}\middle/\int_{\bb R}\frac{r\mu(t)}{(s-t)^2+r^2}\right]=f(s)
\end{equation}
almost everywhere with respect to Lebesgue measure.

We are interested in the behavior of the quotient on the left-hand side of \eqref{221115.1} on a set of full $\mu$-measure, that is, also at points where $\sigma$ is concentrated.
This was achieved in the following theorem whose proof is due to Bj\"orn and Christer Bennewitz \cite{BB-CB}.

\begin{thm}\label{fatou}
Suppose $\mu$ is a non-negative measure on $\bb R$ such that $\int_{\bb R} \mu(t)/(t^2+1)$ is finite.
If $f\in L^\infty(\mu)$, then
$$\lim_{r\downarrow0}\left[\int_{\bb R} \frac{r f(t)\mu(t)}{(s-t)^2+r^2}\middle/ \int_{\bb R} \frac{r\mu(t)}{(s-t)^2+r^2}\right]=f(s)$$
for $\mu$-almost every $s\in\bb R$.
\end{thm}

\begin{proof}
Suppose that $s$ is a $\mu$-Lebesgue point of $f$, i.e., for every positive $\varepsilon$ there is a positive $\delta$ such that
\begin{equation}\label{e:lp}
\int_I |f(t)-f(s)|\mu(t) < \varepsilon \mu(I)
\end{equation}
for any interval $I\subset[s-\delta,s+\delta]$ containing $s$.
Without loss of generality one may assume that $\delta\leq 1$ and that $s-\delta$ and $s+\delta$ are points of continuity for $\mu$.
The differentiation theorem (see, e.g., Theorem B.8.8 of \cite{MR4199125}) implies that $\mu$-almost every point is a $\mu$-Lebesgue point.
Moreover assume that
\begin{equation}\label{e:hmt}
\inf\{\mu((s-r,s+r))/(2r):r>0\}>0.
\end{equation}
Again, this is true for $\mu$-almost all $s$ by Hardy's maximal theorem (see, e.g., Theorem B.8.9 of \cite{MR4199125}) applied to the function $s\mapsto \sup\{2r/\mu((s-r,s+r)): r>0\}$.

Let $F$ be the left-continuous antiderivative of $|f(\cdot)-f(s)|\mu$ satisfying $F(s)=0$ and, similarly, $M$ the left-continuous antiderivative of $\mu$ satisfying $M(s)=0$.
Also abbreviate the expression $r/((s-t)^2+r^2)$ by $p(t)$.
Then the claim may be written as
$$\left.\int_{\bb R} p\,dF\middle/ \int_{\bb R} p\,dM\right.\to 0$$
as $r\downarrow0$.

Equation \eqref{e:lp} may be rephrased in terms of $F$ and $M$ as follows.
If $t>s$
\begin{equation}\label{est1}
F(t)\leq\varepsilon\mu([s,t))=\varepsilon M(t)
\end{equation}
and, if $t<s$,
\begin{equation}\label{est2}
-F(t)\leq\varepsilon\mu([t,s])=\varepsilon(\mu(\{s\})-M(t)).
\end{equation}

Split the integral $\int_{\bb R} p\,dF$ into two parts, namely $I_1=\int_{(s-\delta,s+\delta)} p\,dF$ and $I_2=\int_{(s-\delta,s+\delta)^c} p\,dF$.
Since $p$ is continuous integration by parts yields
$$I_1=(pF)(s+\delta)-(pF)(s-\delta)- \int_{(s-\delta,s+\delta)} F\,dp.$$
Now use the estimates \eqref{est1} and \eqref{est2} to obtain
$$I_1\leq\varepsilon \Big((pM)(s+\delta)-(pM)(s-\delta) +\mu(\{s\})p(s)-\int_{(s-\delta,s+\delta)} M\,dp\Big).$$
Integrating by parts once more one gets
$$I_1\leq \varepsilon \mu(\{s\})p(s)+\varepsilon\int_{(s-\delta,s+\delta)} p\,dM\leq 2\varepsilon \int_{\bb R}p\,dM.$$

To estimate $I_2$ note that, if $|s-t|\geq \delta$ and $\delta\in(0,1]$, then
$$\frac{r(t^2+1)}{(s-t)^2+r^2}\leq 8\frac{s^2+1}{\delta^2}r.$$
To obtain this estimate use $t^2+1\leq 4s^2+2(s^2+1)+2\leq 8(s^2+1)$ when $|t|\leq 2|s|+1$ and $t^2+1\leq 2t^2$ and $(s-t)^2\leq t^2/4$ when $|t|\geq 2|s|+1$.
Therefore
$$I_2\leq 16 \|f\|_\infty\frac{s^2+1}{\delta^2}r\int\frac{\mu(t)}{t^2+1}$$
which tends to $0$ as $r\downarrow0$.
Also
$$\int_{\bb R} p\,dM \geq \int_{(s-r,s+r)} p\,dM \geq \frac{\mu((s-r,s+r))}{2r}$$
which is bounded away from $0$ according to equation \eqref{e:hmt}.
Hence $I_2/\int_{\bb R} p\,dM<\varepsilon$ for sufficiently small $r$.

Combining this estimate with the one for $I_1$ proves the claim, since $\varepsilon$ may be arbitrarily small.
\end{proof}

\bibliographystyle{plain}


\begin{thebibliography}{10}

\bibitem{BB-CB}
{Bj\"orn} Bennewitz and Christer Bennewitz.
\newblock Private communication, 2022.

\bibitem{Bennewitz-sths}
Christer Bennewitz.
\newblock Spectral theory in {H}ilbert space.
\newblock Lecture Notes, 2008.

\bibitem{MR4199125}
Christer Bennewitz, Malcolm Brown, and Rudi Weikard.
\newblock {\em Spectral and scattering theory for ordinary differential
  equations. {V}ol. {I}}.
\newblock Universitext. Springer, Cham, 2020.
\newblock Sturm--Liouville equations.

\bibitem{MR4298818}
Kevin Campbell, Minh Nguyen, and Rudi Weikard.
\newblock On the spectral theory for first-order systems without the unique
  continuation property.
\newblock {\em Linear Multilinear Algebra}, 69(12):2315--2323, 2021.
\newblock Published online: 04 Oct 2019.

\bibitem{MR3046408}
Jonathan Eckhardt, Fritz Gesztesy, Roger Nichols, and Gerald Teschl.
\newblock Weyl-{T}itchmarsh theory for {S}turm-{L}iouville operators with
  distributional potentials.
\newblock {\em Opuscula Math.}, 33(3):467--563, 2013.

\bibitem{MR3095152}
Jonathan Eckhardt and Gerald Teschl.
\newblock Sturm-{L}iouville operators with measure-valued coefficients.
\newblock {\em J. Anal. Math.}, 120:151--224, 2013.

\bibitem{MR1555035}
P.~Fatou.
\newblock S\'{e}ries trigonom\'{e}triques et s\'{e}ries de {T}aylor.
\newblock {\em Acta Math.}, 30(1):335--400, 1906.

\bibitem{MR0068082}
William Feller.
\newblock On second order differential operators.
\newblock {\em Ann. of Math. (2)}, 61:90--105, 1955.

\bibitem{MR914699}
F.~Gesztesy and H.~Holden.
\newblock A new class of solvable models in quantum mechanics describing point
  interactions on the line.
\newblock {\em J. Phys. A}, 20(15):5157--5177, 1987.

\bibitem{MR4047968}
Ahmed Ghatasheh and Rudi Weikard.
\newblock Spectral theory for systems of ordinary differential equations with
  distributional coefficients.
\newblock {\em J. Differential Equations}, 268(6):2752--2801, 2020.

\bibitem{MR0080835}
I.~S. Kac.
\newblock On the existence of spectral functions of certain second-order
  singular differential systems.
\newblock {\em Dokl. Akad. Nauk SSSR (N.S.)}, 106:15--18, 1956.

\bibitem{MR0054078}
M.~G. Kre{\u\i}n.
\newblock On a generalization of investigations of {S}tieltjes.
\newblock {\em Doklady Akad. Nauk SSSR (N.S.)}, 87:881--884, 1952.

\bibitem{MR1397901}
P.~Kurasov.
\newblock Distribution theory for discontinuous test functions and differential
  operators with generalized coefficients.
\newblock {\em J. Math. Anal. Appl.}, 201(1):297--323, 1996.

\bibitem{MR1443392}
P.~Kurasov and J.~Boman.
\newblock Finite rank singular perturbations and distributions with
  discontinuous test functions.
\newblock {\em Proc. Amer. Math. Soc.}, 126(6):1673--1683, 1998.

\bibitem{MR706255}
Angelo~B. Mingarelli.
\newblock {\em Volterra-{S}tieltjes integral equations and generalized ordinary
  differential expressions}, volume 989 of {\em Lecture Notes in Mathematics}.
\newblock Springer-Verlag, Berlin, 1983.

\bibitem{MR4431055}
Steven Redolfi and Rudi Weikard.
\newblock Green's functions for first-order systems of ordinary differential
  equations without the unique continuation property.
\newblock {\em Integral Equations Operator Theory}, 94(2):Paper No. 23, 19,
  2022.

\bibitem{MR1307384}
Marvin Rosenblum and James Rovnyak.
\newblock {\em Topics in {H}ardy classes and univalent functions}.
\newblock Birkh\"{a}user Advanced Texts: Basler Lehrb\"{u}cher. [Birkh\"{a}user
  Advanced Texts: Basel Textbooks]. Birkh\"{a}user Verlag, Basel, 1994.

\bibitem{MR1756602}
A.~M. Savchuk and A.~A. Shkalikov.
\newblock Sturm-{L}iouville operators with singular potentials.
\newblock {\em Mathematical Notes}, 66(6):741--753, 1999.
\newblock Translated from Mat. Zametki, Vol. 66, pp. 897--912 (1999).

\end{thebibliography}
\end{document}